\documentclass[12pt]{article}
\usepackage{amssymb}
\usepackage{amsmath}
\usepackage{amsthm}
\usepackage{color}
\usepackage{comment}
\usepackage{geometry}		
\usepackage{graphicx}

\makeatletter
	\@addtoreset{equation}{section}
\makeatother
\renewcommand{\labelenumi}{\roman{enumi})}

\let\originalleft\left
\let\originalright\right
\renewcommand{\left}{\mathopen{}\mathclose\bgroup\originalleft}
\renewcommand{\right}{\aftergroup\egroup\originalright}

\begin{document}

\newcommand\cX{\mathcal{X}}
\newcommand{\rD}{{\rm D}}
\newcommand{\re}{{\rm e}}

\newcommand{\removableFootnote}[1]{\footnote{#1}}

\newtheorem{theorem}{Theorem}[section]
\newtheorem{corollary}[theorem]{Corollary}
\newtheorem{lemma}[theorem]{Lemma}
\newtheorem{proposition}[theorem]{Proposition}

\theoremstyle{definition}
\newtheorem{definition}{Definition}[section]
\newtheorem{example}[definition]{Example}

\theoremstyle{remark}
\newtheorem{remark}{Remark}[section]



\title{The two-dimensional border-collision normal form with a zero determinant}
\author{
D.J.W.~Simpson\\\\
School of Mathematical and Computational Sciences\\
Massey University\\
Palmerston North, 4410\\
New Zealand
}
\maketitle


\begin{abstract}

The border-collision normal form is a piecewise-linear family of continuous maps
that describe the dynamics near border-collision bifurcations.
Most prior studies assume each piece of the normal form is invertible,
as is generic from an abstract viewpoint,
but in applied problems one piece of the map often has degenerate range,
corresponding to a zero determinant.
This provides simplification, yet even in two dimensions the dynamics can be incredibly rich.
The purpose of this paper is to determine broadly how the dynamics of
the two-dimensional border-collision normal form with a zero determinant
differs for different values of its parameters.
We identify parameter regions of period-adding,
period-incrementing, mode-locking, and component doubling of chaotic attractors,
and characterise the dominant bifurcation boundaries.
The intention is for the results to enable border-collision bifurcations
in mathematical models to be analysed more easily and effectively,
and we illustrate this with a flu epidemic model and two stick-slip friction oscillator models.
We also describe three novel bifurcation structures that remain to be explored.

\end{abstract}

\section{Introduction}
\label{sec:intro}

To understand the behaviour of a dynamical system across a range of parameter values,
it is essential to identify critical parameter values --- bifurcations --- where the dynamics changes in a fundamental way.
Classical bifurcation theory relies on the equations of motion being smooth \cite{Ku04}.
But mathematical models often involve nonsmooth equations,
e.g.~to describe impacts in mechanical systems \cite{Br99,WiDe00},
switches in control systems \cite{Li03},
and decisions in ecology, economics, and society \cite{DeGr07,PuSu06}.
Such models exhibit novel bifurcations, most simply border-collision bifurcations (BCBs)
whereby a fixed point of a piecewise-smooth map collides with a switching manifold.
BCBs include grazing events of limit cycles
and represent the onset of chaos in diverse applications \cite{DiBu08}.

\begin{figure}[b!]
\begin{center}
\includegraphics[height=7cm]{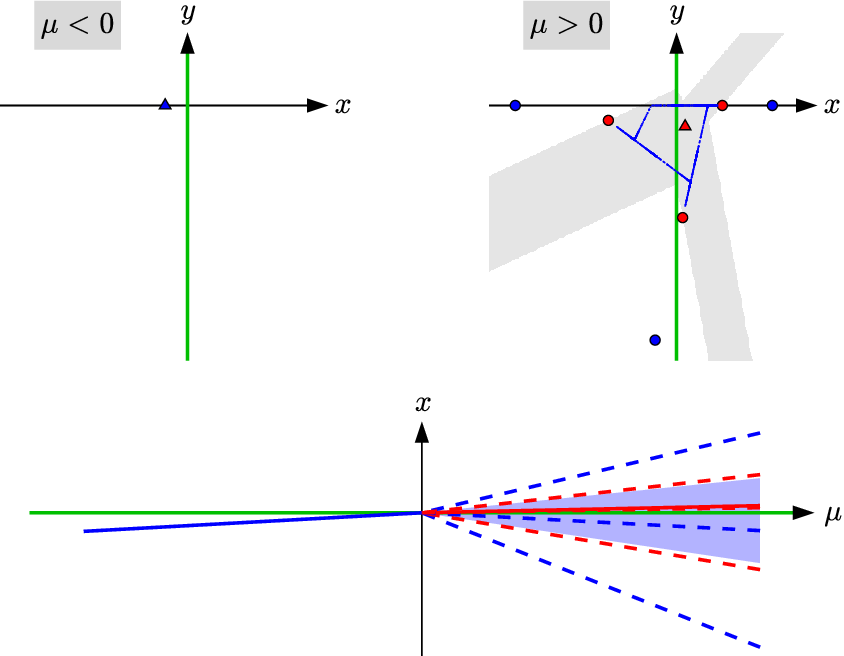}
\caption{
The upper plots are phase portraits of
the two-dimensional border-collision normal form \eqref{eq:bcnf}
with $(\tau_L,\delta_L,\tau_R,\delta_R) = (-0.4,0,-0.55,2.1)$.
The map has a stable fixed point (blue triangle) for $\mu < 0$ and an unstable fixed point (red triangle) for $\mu > 0$.
For $\mu > 0$ it also has a stable period-three solution (blue circles) and a chaotic attractor (blue dots).
The basin of attraction (grey) of the chaotic attractor
is bounded by the stable manifold of a saddle period-three solution (red circles).
The lower plot is a bifurcation diagram showing how the $x$-values of the invariant sets
vary with $\mu$ (solid lines: fixed points;
dashed lines: period-three solutions;
shaded area: chaotic attractor).
\label{fig:zRbifDiagSchem}
} 
\end{center}
\end{figure}

This paper concerns BCBs for maps that are continuous and piecewise-linear, to leading order.
By truncating the map and changing coordinates we obtain the border-collision normal form \cite{Di03,Si16}.
In two dimensions this form can be written as
\begin{equation}
\begin{bmatrix} x \\ y \end{bmatrix} \mapsto
\begin{cases}
\begin{bmatrix} \tau_L x + y + \mu \\ -\delta_L x \end{bmatrix}, & x \le 0, \\
\begin{bmatrix} \tau_R x + y + \mu \\ -\delta_R x \end{bmatrix}, & x \ge 0,
\end{cases}
\label{eq:bcnf}
\end{equation}
where $(x,y) \in \mathbb{R}^2$ is the system state,
$\mu \in \mathbb{R}$ is the primary bifurcation parameter,
and $\tau_L, \delta_L, \tau_R, \delta_R \in \mathbb{R}$ are additional parameters.
The line $x=0$ is the switching manifold,
and the normal form exhibits a border-collision bifurcation by varying $\mu$ through zero.
Fig.~\ref{fig:zRbifDiagSchem} shows a typical example.
Here the map has a stable fixed point for all $\mu < 0$,
and a stable period-three solution and a chaotic attractor for all $\mu > 0$.
Since \eqref{eq:bcnf} is piecewise-linear,
any bounded invariant set of \eqref{eq:bcnf} collapses linearly to $(x,y) = (0,0)$ as $\mu$ tends monotonically to $0$.

Different dynamics occurs for different values of $\tau_L$, $\delta_L$, $\tau_R$, and $\delta_R$.
These values are the traces and determinants of the
Jacobian matrix of \eqref{eq:bcnf} on each side of $x = 0$.
Mostly we are interested in attracting invariant sets
as these dictate the long-term behaviour of typical orbits.

To analyse a BCB in a mathematical model, the normal form \eqref{eq:bcnf} can be utilised as follows.
We first determine the values $\tau_L$, $\delta_L$, $\tau_R$, and $\delta_R$.
To do this it is not necessary to explicitly truncate and perform the change of coordinates.
We simply evaluate the Jacobian matrix of each piece of the map at the bifurcation,
say by using finite differences,
then set $\tau_L$, $\delta_L$, $\tau_R$, and $\delta_R$ equal to the traces and determinants of these matrices.
The idea is to next invoke known theory for the dynamics of \eqref{eq:bcnf},
then claim that these dynamics occur in the model.
Indeed the terms omitted by truncation typically have no qualitative effect on the local dynamics.
By working with \eqref{eq:bcnf}, codimension-two points on curves of BCBs
can be identified accurately because the truncation removes nearby bifurcations and features
that may obscure the dynamics specific to the BCB.

This approach relies on a good prior understanding of the 
dynamics of \eqref{eq:bcnf} for all $\tau_L$, $\delta_L$, $\tau_R$, and $\delta_R$.
This has motivated many studies, firstly that of Nusse and Yorke \cite{NuYo92}
where \eqref{eq:bcnf} was introduced
and where it was shown that BCBs
can create a plethora of dynamical transitions that are not possible for smooth systems.
Later Banerjee and Grebogi \cite{BaGr99} focussed on the
dissipative, orientation-preserving case $0 < \delta_L, \delta_R < 1$,
and found eleven basic qualitative types of transitions.
Around this time it became clear that a complete characterisation of the dynamics of \eqref{eq:bcnf}
is almost certainly unattainable; there are too many possibilities.
Subsequent studies concentrated on particular features and patterns,
such as Arnold tongues, which have a distinctive sausage-string structure \cite{Si17c,Si24d,SuGa08,ZhMo06b},
and chaotic attractors, which persist throughout open regions of parameter space \cite{BaYo98,GhMc23,Gl17}.
Overviews of the bifurcation structure have recently been presented
for the non-invertible case $\delta_L \delta_R < 0$ \cite{FaSi23},
and the constant determinant case $\delta_L = \delta_R$ \cite{SuAv22}.

But in many applications BCBs have $\delta_L = 0$ or $\delta_R = 0$.
This occurs for grazing-sliding bifurcations,
where a limit cycle of a Filippov system grazes a discontinuity surface \cite{DiKo02},
and event collisions for systems with delayed switching,
where the time taken for a limit cycle to cross and then return to a switching surface
matches the delay time \cite{SiKo10b}.
A determination of the dynamics and bifurcations of \eqref{eq:bcnf} with $\delta_L = 0$ or $\delta_R = 0$
is needed, and such is the aim of this paper.

Previous studies of \eqref{eq:bcnf} with $\delta_L = 0$ or $\delta_R = 0$ have treated isolated aspects of the dynamics.
Kowalczyk \cite{Ko05} computed fixed points and period-two solutions
and identified parameter regions where the map has an attractor contained in the union of two or three line segments.
He argued the attractor is chaotic by showing no stable periodic solutions exist,
and observed that the number of connected components that comprise the attractor differs across the parameter regions.
Later Szalai and Osinga \cite{SzOs08} analysed invariant polygons in a two-parameter subfamily.
They used symbolic dynamics to argue that if the attractor leaks out of the polygon it can be chaotic,
and in \cite{SzOs09} studied Arnold tongues having the sausage-string geometry.

\subsection{Outline}

This paper starts in \S\ref{sec:skewTentMaps} with
the dynamics of one-dimensional piecewise-linear maps (skew tent maps)
as such maps result from setting {\em both} $\delta_L = 0$ and $\delta_R = 0$ in the normal form.
It is instructive to see how the additional spatial dimension realised by varying one of $\delta_L$ and $\delta_R$ from zero
creates new possibilities, such as coexisting attractors, while some key features
of the one-dimensional setting are retained, such as component doubling.

In \S\ref{sec:grazingSliding} we review grazing-sliding bifurcations.
We show how the oscillatory dynamics local to grazing-sliding bifurcations of three-dimensional Filippov systems
are described by \eqref{eq:bcnf} with $\delta_L = 0$ or $\delta_R = 0$.

In the remainder of the paper we fix $\delta_L = 0$, without loss of generality (see Remark \ref{rm:LRswap}).
In this case \eqref{eq:bcnf} reduces to
\begin{equation}
\begin{bmatrix} x \\ y \end{bmatrix} \mapsto
\begin{cases}
\begin{bmatrix} \tau_L x + y + \mu \\ 0 \end{bmatrix}, & x \le 0, \\
\begin{bmatrix} \tau_R x + y + \mu \\ -\delta_R x \end{bmatrix}, & x \ge 0,
\end{cases}
\label{eq:f}
\end{equation}
and in \S\ref{sec:basics} we describe elementary aspects of \eqref{eq:f}.

Sections \ref{sec:neg} and \ref{sec:pos} concern \eqref{eq:f} with $\mu < 0$ and $\mu > 0$ respectively.
We combine a comprehensive numerical analysis
with explicit calculations of the dominant bifurcations
and characterise the attractors of \eqref{eq:f} across parameter space.
We derive formulas for the boundaries of periodicity regions
and some homoclinic and heteroclinic bifurcations that serve as crises 
where a chaotic attractor jumps in size or is destroyed.
We use renormalisation to identify bifurcation surfaces where a chaotic attractor undergo component doubling.
Chaotic attractors are in many places created along curves of {\em shrinking points}.
Several of these admit closed-form expressions
which is not the case when $\delta_L \ne 0$ and $\delta_R \ne 0$ \cite{SiMe08b}.
We also identify three periodicity structures
residing in seas of otherwise chaotic dynamics that have yet to be explored.
To the authors' knowledge two of these have not been reported before.

In \S\ref{sec:appl} we apply the results to a flu epidemic model of Roberts {\em et al.}~\cite{RoHi19b}
and two stick-slip friction models.
For each we compute the corresponding values of $\tau_L$, $\tau_R$, and $\delta_R$
and use the results of \S\ref{sec:neg} and \S\ref{sec:pos} explain the model dynamics.
This shows how the dynamics described in previous works fits into a bigger picture.
Other models to which the results could be applied
include the ecological model of Zhou and Tang \cite{ZhTa22}.
Summarising remarks are provided in \S\ref{sec:conc}.

\section{Skew tent maps}
\label{sec:skewTentMaps}

In this section we study the skew tent map family
\begin{equation}
x \mapsto \begin{cases}
s_L x + \eta, & x \le 0, \\
s_R x + \eta, & x \ge 0,
\end{cases}
\label{eq:skewTentMap}
\end{equation}
where $s_L, s_R, \eta \in \mathbb{R}$ are parameters.
As with $\mu$ in the normal form \eqref{eq:bcnf},
the structure of the dynamics of \eqref{eq:skewTentMap} depends only on the sign of $\eta$. 

\begin{figure}[b!]
\begin{center}
\includegraphics[width=15.6cm]{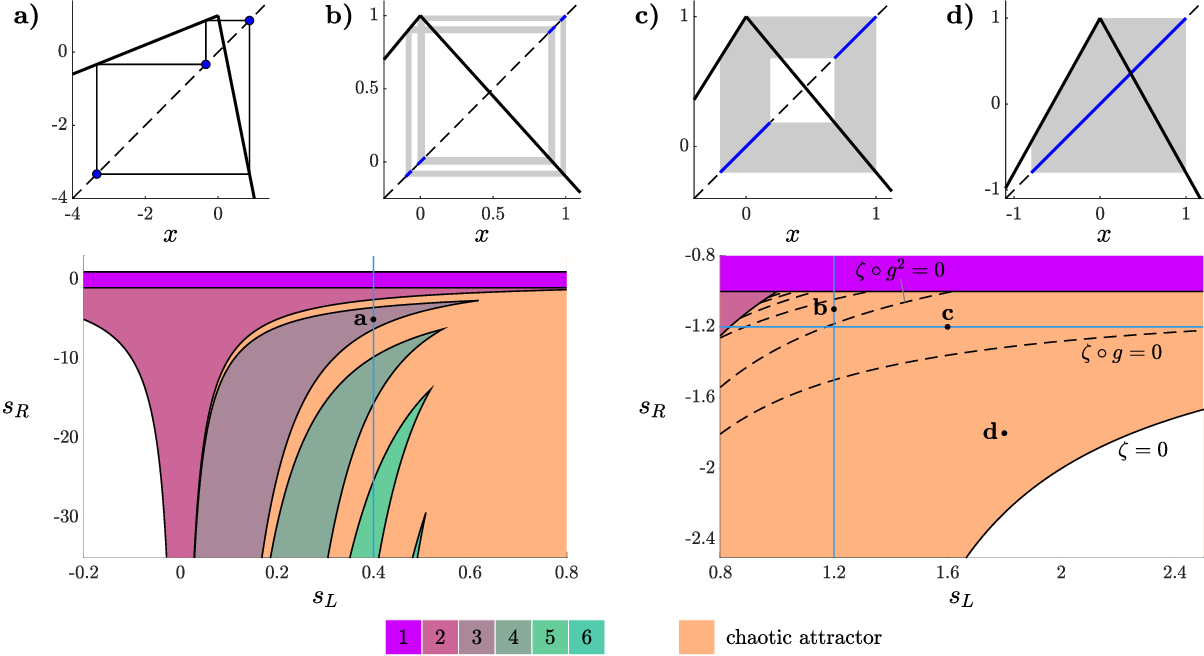}
\caption{
Cobweb diagrams and two-parameter bifurcation diagrams
of the skew tent map family \eqref{eq:skewTentMap} with $\eta = 1$.
The orange regions are where the map has a chaotic attractor.
The other coloured regions are where the map has a stable period-$p$ solution with exactly one point in $x > 0$
for the value of $p$ indicated by the colour bar.
In the white regions the map has no attractor.
The cobweb diagrams show the attractor and use:
(a) $(s_L,s_R) = (0.4,-5)$;
(b) $(s_L,s_R) = (1.2,-1.1)$;
(c) $(s_L,s_R) = (1.6,-1.2)$;
(d) $(s_L,s_R) = (1.8,-1.8)$.
The blue lines indicate parameter points arising as special cases in later sections.
\label{fig:zRbifSetSkewTentMap}
} 
\end{center}
\end{figure}

Fig.~\ref{fig:zRbifSetSkewTentMap} summarises the attractor of \eqref{eq:skewTentMap} for $\eta > 0$.
The lower plots show how the nature of the attractor differs across the $(s_L,s_R)$-plane.
The upper plots are sample cobweb diagrams in which the attractor is indicated by blue points on the $45^\circ$-degree line.
For details, including formulas for the curves in the lower plots, 
refer to \cite{ItTa79b,MaMa93,NuYo95,SuAv16}.
The case $\eta < 0$ can be accommodated via the substitution $x \mapsto -x$ that effectively swaps $s_L$ and $s_R$.

We now highlight three features of \eqref{eq:skewTentMap} that will be important in later sections.

\begin{remark}
If \eqref{eq:skewTentMap} has an attractor this attractor is unique.
This can be inferred from an argument that any basin of attraction must contain a neighbourhood of $x=0$,
or from the fact that \eqref{eq:skewTentMap} does not have positive Schwarzian derivative at any points \cite{MiTh88}.
\end{remark}

\begin{remark}
If \eqref{eq:skewTentMap} has a stable periodic solution,
this solution has exactly one point in $x > 0$ if $\eta > 0$,
and exactly one point in $x < 0$ if $\eta < 0$.
For example, with $\eta = 1$ and $(s_L,s_R) = (0.4,-5)$, as in Fig.~\ref{fig:zRbifSetSkewTentMap}a,
the map has a stable period-three solution with one point in $x > 0$.
The lower-left plot in Fig.~\ref{fig:zRbifSetSkewTentMap}
shows regions where \eqref{eq:skewTentMap} has stable periodic solutions up to period six.
All higher periods occur at more negative values of $s_R$.
\end{remark}

\begin{remark}
As we move about the orange region where \eqref{eq:skewTentMap} has a chaotic attractor,
the number of intervals that comprise the attractor doubles (or halves) as we cross the dashed curves.
These curves can be computed by renormalisation as follows.
The curve $\zeta(\tau_L,\tau_R) = 0$, where
\begin{equation}
\zeta(s_L,s_R) = s_L s_R + s_L - s_R \,,
\label{eq:divergenceFunction}
\end{equation}
is where the second iterate of the critical point $x=0$ coincides with the fixed point $\frac{\eta}{1-s_L}$ of the left piece of the map.
Below this curve \eqref{eq:skewTentMap} has no attractor.
Immediately above this curve \eqref{eq:skewTentMap} has a single interval attractor, as in Fig.~\ref{fig:zRbifSetSkewTentMap}d.
As explained in \cite{GhMc24,ItTa79b,VeGl90},
as we cross the curve $\zeta \left( g^k(s_L,s_R) \right) = 0$, where $k \ge 1$ and $g$ is the {\em renormalisation operator}
\begin{equation}
g(s_L,s_R) = \left( s_R^2, s_L s_R \right),
\label{eq:renormalisationOperator}
\end{equation}
the number of intervals comprising the attractor changes from $2^{k-1}$ to $2^k$.
For $k = 1,\ldots,5$ these curves are shown dashed in the lower-right plot of Fig.~\ref{fig:zRbifSetSkewTentMap}.
\end{remark}

\section{Grazing-sliding bifurcations}
\label{sec:grazingSliding}

In this section we explain how the normal form with a zero determinant arises
for grazing-sliding bifurcations in Filippov systems.
This was shown originally by di Bernardo {\em et al.}~\cite{DiKo02,DiKo03},
and for further details refer to the textbook \cite{DiBu08}.

Consider a three-dimensional ODE system of the form
\begin{equation}
\frac{d X}{d t} = \begin{cases}
F_L(X), & H(X) < 0, \\
F_R(X), & H(X) > 0,
\end{cases}
\label{eq:Filippov}
\end{equation}
where $X(t) \in \mathbb{R}^3$ is the system state,
$F_L$ and $F_R$ are smooth vector fields,
and $H : \mathbb{R}^3 \to \mathbb{R}$ is a smooth switching function.
We assume the {\em discontinuity surface}
\begin{equation}
\Sigma = \left\{ X \in \mathbb{R}^3 \,\middle|\, H(X) = 0 \right\}
\nonumber
\end{equation}
is a smooth two-dimensional manifold.
 
\begin{figure}[b!]
\begin{center}
\includegraphics[height=6cm]{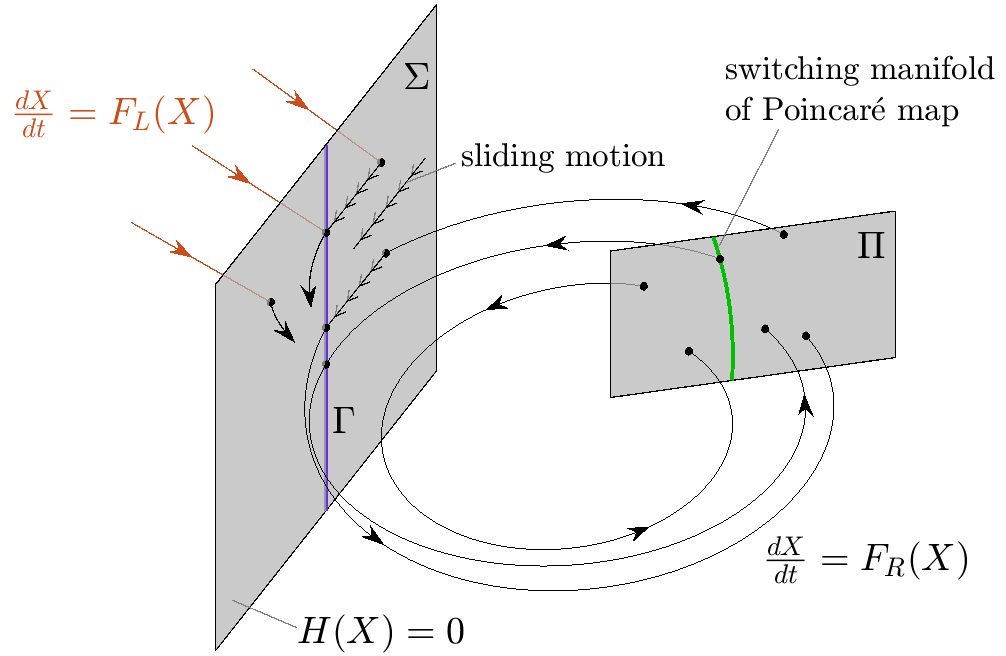}
\caption{
A sketch of the phase space of a three-dimensional Filippov system \eqref{eq:Filippov}
illustrating the construction of the Poincar\'e map $P$.
\label{fig:zRgrazSchem}
} 
\end{center}
\end{figure}

For our purposes we assume $F_L$ directs solutions towards $\Sigma$,
while $F_R$ defines a curve $\Gamma$ of {\em visible folds}
where orbits of $F_R$ graze $\Sigma$ quadratically \cite{Je18b}.
This situation is depicted in Fig.~\ref{fig:zRgrazSchem}.

On one side of $\Gamma$ solutions are directed towards $\Sigma$ by both $F_L$ and $F_R$.
This part of $\Sigma$ is called an {\em attracting sliding region}.
We assume solutions evolve on this region
according to $\frac{d X}{d t} = F_S(X)$, where $F_S$ is a {\em sliding vector field} constructed from $F_L$ and $F_R$.
The standard way to define $F_S$ is as the unique convex combination of $F_L$ and $F_R$ that is tangent to $\Sigma$,
and with this convention \eqref{eq:Filippov} is termed a {\em Filippov system}.
Note that sliding motion ends when solutions reach $\Gamma$.

Now let $\Pi$ be a smooth cross-section of phase space in $H(X) > 0$
that is transverse to $F_R$, see again Fig.~\ref{fig:zRgrazSchem}.
Assume that the forward orbits of all points in some subset $S \subset \Pi$ return to $\Pi$
after motion either entirely within $H(X) > 0$ or with one sliding segment.
This defines a two-dimensional Poincar\'e map $P$ that we can write as
\begin{equation}
P(u) = \begin{cases}
P_L(u), & J(u) \le 0, \\
P_R(u), & J(u) \ge 0.
\end{cases}
\label{eq:PoincareMap}
\end{equation}
Here $P_L$ corresponds to orbits with a sliding segment, while $P_R$ corresponds to orbits without a sliding segment.
The map $P$ is continuous across its switching manifold $J(u) = 0$, which corresponds to orbits that graze $\Sigma$ on $\Gamma$.
It is clear that $P_R$ is $C^1$, since $F_R$ is smooth and intersects $\Pi$ transversely.
It turns out that $P_L$ is also $C^1$ --- this follows from asymptotic calculations
and relies on the fact that $F_S = F_R$ on $\Gamma$, see \cite{DiKo02} --- hence
the Poincar\'e map does not have square-root singularity,
as is the case for some other grazing scenarios \cite{DiBu01,FrNo97,No91}.
Importantly, the range of $P_L$ is one-dimensional
because all orbits with sliding segments return to $\Pi$ by passing through $\Gamma$.

A {\em grazing-sliding bifurcation} occurs when a limit cycle of $F_R$ hits $\Gamma$ as system parameters are varied.
In the context of the Poincar\'e map, this corresponds to a fixed point of $P_R$ hitting the switching manifold $J(u) = 0$.
By replacing $P_L$ and $P_R$ with their linearisations evaluated at the bifurcation point,
we produce a continuous piecewise-linear map whose dynamics approximates the local dynamics of $P$.
The piecewise-linear map can be converted to the normal form \eqref{eq:bcnf}, assuming an observability condition is satisfied \cite{Si16},
and preserving the notions of left and right.
Notice $\delta_L = 0$ because $P_L$ has degenerate range, hence we have \eqref{eq:f}.

\begin{remark}
Suppose instead $\delta_R = 0$ in the normal form \eqref{eq:bcnf}.
Under the substitution $(x,y) \mapsto (-x,-y)$
we obtain \eqref{eq:f} except with $(\mu,\tau_L,\tau_R,\delta_R)$
in place of $(-\mu,\tau_R,\tau_L,\delta_L)$.
So for example the dynamics of \eqref{eq:bcnf} with $(\tau_L,\delta_L,\tau_R,\delta_R) = (1,2,3,0)$
are the same of those of \eqref{eq:f} with $(\tau_L,\tau_R,\delta_R) = (3,1,2)$,
but occurring for opposite signs of $\mu$.
In this way our results in later sections
can be applied to the normal form \eqref{eq:bcnf} with $\delta_R = 0$
by switching $L$ and $R$ and the sign of $\mu$.
\label{rm:LRswap}
\end{remark}

\section{Fixed points and symbolic representations of periodic solutions}
\label{sec:basics}

Recall that in applied settings, e.g.~grazing-sliding bifurcations,
\eqref{eq:f} approximates the dynamics of a more general piecewise-smooth system.
In this context the approximation is only reasonable for small values of $x$, $y$, and $\mu$.
But since \eqref{eq:f} is piecewise-linear,
for fixed $\tau_L$, $\tau_R$, and $\delta_R$
the structure of its dynamics depends only on the sign of $\mu$, see again Fig.~\ref{fig:zRbifDiagSchem}.
So if $\mu < 0$ it suffices to consider $\mu = -1$, treated in \S\ref{sec:neg},
while if $\mu > 0$ it suffices to consider $\mu = 1$, treated in \S\ref{sec:pos}.

Let
\begin{equation}
f_L(x,y;\mu) = \begin{bmatrix} \tau_L x + y + \mu \\ 0 \end{bmatrix}, \qquad
f_R(x,y;\mu) = \begin{bmatrix} \tau_R x + y + \mu \\ -\delta_R x \end{bmatrix},
\nonumber
\end{equation}
denote the pieces of \eqref{eq:f}.
Also let
\begin{equation}
A_L = \rD f_L(x,y;\mu) = \begin{bmatrix} \tau_L & 1 \\ 0 & 0 \end{bmatrix}, \qquad
A_R = \rD f_R(x,y;\mu) = \begin{bmatrix} \tau_R & 1 \\ -\delta_R & 0 \end{bmatrix},
\nonumber
\end{equation}
denote the left and right Jacobian matrices.
For any $\mu \in \mathbb{R}$ the image of the left half-plane $x \le 0$ under $f_L$ is the $x$-axis,
while the image of the right half-plane $x \ge 0$ under $f_R$ is the lower half-plane if $\delta_R > 0$,
and is the upper half-plane if $\delta_R < 0$.
This explains why in Fig.~\ref{fig:zRbifDiagSchem}, which has $\delta_R > 0$,
all invariant sets are contained in the lower half-plane
and include points on the $x$-axis.

\subsection{Fixed points}

If $\tau_R \ne 1 + \delta_R$ then $f_R$ has the unique fixed point
\begin{equation}
\begin{bmatrix} x^R \\ y^R \end{bmatrix} = \frac{\mu}{1 - \tau_R + \delta_R} \begin{bmatrix} 1 \\ -\delta_R \end{bmatrix}.
\label{eq:zR}
\end{equation}
If $x^R > 0$ then $(x^R,y^R)$ is a fixed point of \eqref{eq:f} and said to be {\em admissible};
if $x^R < 0$ then $(x^R,y^R)$ is {\em virtual}.
If $(x^R,y^R)$ is admissible then it is asymptotically stable if and only if
both eigenvalues of $\rD f_R(x^R,y^R) = A_R$ have modulus less than $1$.
Since $\tau_R$ and $\delta_R$ are the trace and determinant of $A_R$,
this occurs when $|\tau_R| - 1 < \delta_R < 1$.

The same calculations apply to $f_L$, but are simpler because $\delta_L = 0$.
If $\tau_L \ne 1$ then $f_L$ has the unique fixed point
\begin{equation}
\begin{bmatrix} x^L \\ y^L \end{bmatrix} = \frac{\mu}{1 - \tau_L} \begin{bmatrix} 1 \\ 0 \end{bmatrix}.
\label{eq:zL}
\end{equation}
If $x^L < 0$ then $(x^L,y^L)$ is admissible; if $x^L > 0$ then $(x^L,y^L)$ is virtual.
If $(x^L,y^L)$ is admissible then it is asymptotically stable if and only if $|\tau_L| < 1$.
These observations can be summarised as follows.

\begin{proposition}
If $\tau_R \ne 1 + \delta_R$
then $(x^R,y^R)$ is admissible and asymptotically stable
if and only if $|\tau_R| - 1 < \delta_R < 1$ and $\mu > 0$.
If $\tau_L \ne 1$
then $(x^L,y^L)$ is admissible and asymptotically stable
if and only if $|\tau_L| < 1$ and $\mu < 0$.
\end{proposition}

With $\mu = 0$ the origin $(0,0)$ is a fixed point of \eqref{eq:f}.
The asymptotic stability of this point is difficult to characterise in terms of $\tau_L$, $\tau_R$, and $\delta_R$
because the map is not differentiable at this point \cite{Si20d}.
However, the dynamics of \eqref{eq:f} with $\mu = 0$ will not concern us as it constitutes a special case.

\subsection{Periodic solutions and symbolic representations}

We now explain how periodic solutions of \eqref{eq:f} can be encoded symbolically
using {\em words} comprised of $L$'s and $R$'s.
We cover this material fairly briefly; for a more detailed exposition refer to \cite{Si16,SiMe09}.

It is convenient to denote the map \eqref{eq:f} by $f(x,y)$.
Suppose $f$ has a period-$p$ solution $\{ (x_0,y_0), (x_1,y_1), \ldots, (x_{p-1},y_{p-1}) \}$,
where
\begin{equation}
f(x_0,y_0) = (x_1,y_1),~
f(x_1,y_1) = (x_2,y_2),~
\ldots,~
f(x_{p-1},y_{p-1}) = (x_0,y_0).
\nonumber
\end{equation}
Suppose the periodic solution has no points on the switching manifold,
and define $\cX = \cX_0 \cX_1 \cdots \cX_{p-1}$
by $\cX_i = L$ if $x_i < 0$, and $\cX_i = R$ if $x_i > 0$, for all $i$.
We then refer to the periodic solution as a {\em $\cX$-cycle}.
Note we often use powers to write words more succinctly,
e.g.~$L R^2$ as shorthand for $LRR$.

Since the $\cX$-cycle has no points on the switching manifold,
$f$ is differentiable in a neighbourhood of the $\cX$-cycle.
Thus it has well-defined stability multipliers
given by the eigenvalues of the Jacobian matrix $\left( \rD f^p \right)(x_0,y_0)$.
This matrix is
\begin{equation}
M_\cX = A_{\cX_{p-1}} \cdots A_{\cX_1} A_{\cX_0} \,,
\label{eq:MX}
\end{equation}
so the $\cX$-cycle is asymptotically stable if and only if both eigenvalues of $M_\cX$ have modulus less than $1$.

In the next two sections we identify parameter regions
where $f$ has an asymptotically stable $\cX$-cycle for various words $\cX$.
Boundaries of these regions can be grouped into three types.
First, there are boundaries where the $\cX$-cycle loses stability
by $M_\cX$ attaining an eigenvalue of $-1$
or a complex conjugate pair of eigenvalues with modulus $1$.
Second, there are boundaries where $M_\cX$ has an eigenvalue $1$.
These boundaries are best interpreted as a loss of existence rather than a loss of stability
because as we approach the boundary all points of the $\cX$-cycle tend to infinity \cite{Si16}.
Thirdly, there are boundaries where one point of the $\cX$-cycle hits the switching manifold.
These are examples of BCBs.

\section{Bifurcation structures for $\mu < 0$}
\label{sec:neg}

\begin{figure}[b!]
\begin{center}
\includegraphics[width=15.6cm]{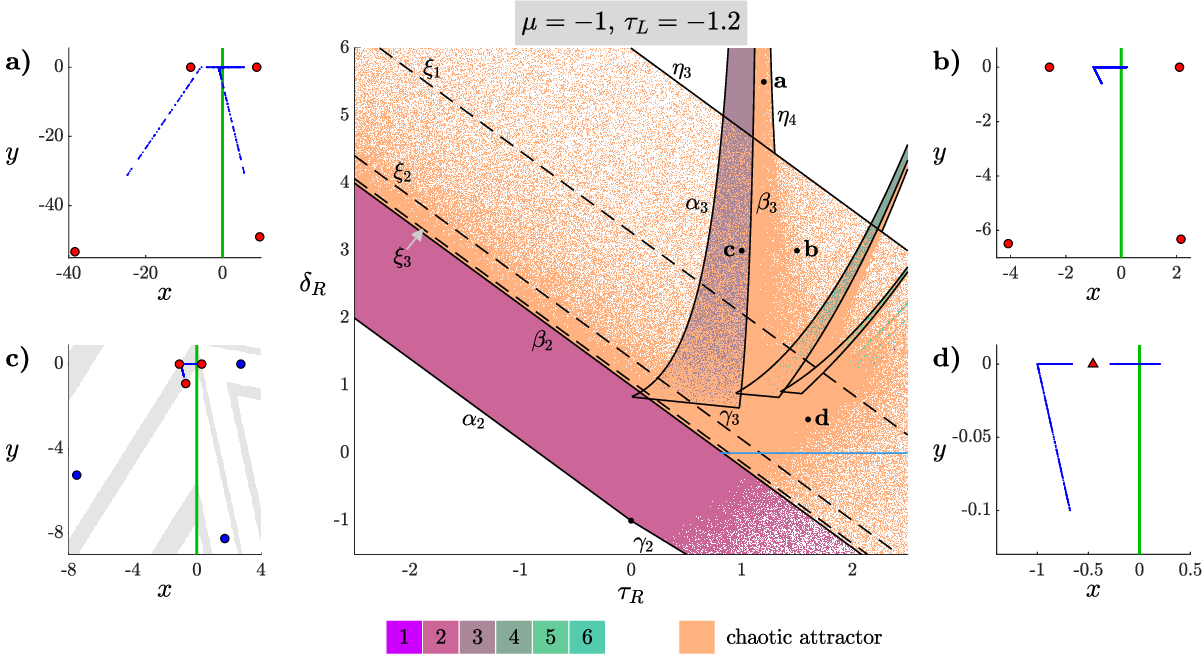}
\caption{
A two-parameter bifurcation diagram of \eqref{eq:f} with $\tau_L = -1.2$ and $\mu = -1$.
The colours indicate the nature of the attractor and are the result of numerics described in \S\ref{sub:numerics}.
The solid black curves are boundaries of regions where there exists an
asymptotically stable $L R^{p-1}$-cycle, for $p = 2,3,4,5$;
also $\eta_3$ is the line \eqref{eq:eta3}.
The dashed lines $\xi_k$ are where $\zeta \left( g^{k-1} \left( g_-(\tau_L,\tau_R,\delta_R) \right) \right) = 0$, for $k = 1,2,3$,
where the number of connected components of the chaotic attractor doubles (or halves).
The blue line at $\delta_R = 0$ corresponds to the
horizontal blue line in the lower-right plot of Fig.~\ref{fig:zRbifSetSkewTentMap}.
The phase portraits correspond to the $(\tau_R,\delta_R)$ points
{\bf a}:~$(1.2,5.5)$;
{\bf b}:~$(1.5,3)$;
{\bf c}:~$(1,3)$;
{\bf d}:~$(1.6,0.5)$.
These show attractors in blue and a selected unstable periodic solution in red.
In panel (c) the basin of attraction of the chaotic attractor is shaded grey.
\label{fig:zRbifSetNegative}
} 
\end{center}
\end{figure}

In this section we describe the dynamics of \eqref{eq:f} with $\mu < 0$.
If $\tau_L > 1$ numerical explorations suggest \eqref{eq:f} cannot have an attractor.
If $-1 < \tau_L < 1$ the fixed point $(x^L,y^L)$ is admissible and asymptotically stable
(although it may not be the only attractor).
Thus we focus on the case $\tau_L < -1$.
Fig.~\ref{fig:zRbifSetNegative} summarises the nature of the attractor when $\tau_L = -1.2$.
Other values of $\tau_L < -1$ yield qualitatively similar pictures.
We now work towards explaining Fig.~\ref{fig:zRbifSetNegative}.

\subsection{Numerical methods}
\label{sub:numerics}

We first describe how the two-parameter bifurcation diagram in the middle of Fig.~\ref{fig:zRbifSetNegative} was computed
(later figures use the same methodology).
For each point in a $1000 \times 1000$ grid of equi-spaced $(\tau_R,\delta_R)$ points,
we numerically computed $M = 10^5$ iterates of the forward orbit of a random initial point $(x_0,y_0)$.
If the last computed iterate $(x_M,y_M)$ had Euclidean norm exceeding $10^5$,
we assumed the orbit was diverging and coloured the $(\tau_R,\delta_R)$ point white.
Otherwise we iterated $30$ more times to produce points $(x_{M+i},y_{M+i})$ for $i = 1,2,\ldots,30$.
If there was $i$ such that $(x_{M+i},y_{M+i}) - (x_M,y_M)$ had Euclidean norm less than $10^{-10}$,
we concluded the orbit was converging to a period-$p$ solution, where $p$ is the smallest such value of $i$,
and coloured the $(\tau_R,\delta_R)$ point according to the colour bar.
If no such $i$ was detected we coloured the point orange if
a numerically computed maximal Lyapunov exponent exceeded $0.001$,
suggesting chaos, and yellow otherwise,
suggesting the orbit is converging to either a quasi-periodic solution
or a periodic solution with period greater than $30$
(although no such points are present in Fig.~\ref{fig:zRbifSetNegative}).

This brute-force numerical approach appears to find attractors relatively effectively.
By using random initial points the output reveals regions
where multiple attractors are present and regions where attractors are not globally attracting.
For example, the middle of Fig.~\ref{fig:zRbifSetNegative} is speckled orange and grey
corresponding to the coexistence of chaotic and period-three attractors, such as at the point (c).
To left of the curve $\alpha_3$, defined below,
the diagram is instead speckled orange and white
corresponding to the presence of a chaotic attractor, but not one that is globally attracting.

The bifurcation diagram has been overlaid with several black curves.
These are bifurcation boundaries plotted by using explicit formulas described below.

\subsection{Periodic solutions, chaos, and crises}
\label{sub:crises}

At parameter point {\bf c}, \eqref{eq:f} has chaotic and period-three attractors.
The period-three attractor is an $L R^2$-cycle (shown as blue circles in Fig.~\ref{fig:zRbifSetNegative}c).
Its basin of attraction is bounded by the stable manifold of a saddle-type $L^2 R$-cycle
(shown as red circles in Fig.~\ref{fig:zRbifSetNegative}c).

At {\bf c} the chaotic attractor lies fairly close to the $L^2 R$-cycle.
Indeed as we move upwards from {\bf c} the chaotic attractor collides with the $L^2 R$-cycle on the line $\eta_3$ and is destroyed.
Here the stable and unstable manifolds of the $L^2 R$-cycle develop non-trivial intersections.
This type of {\em crisis} \cite{GrOt83} is a piecewise-linear analogue of a first homoclinic tangency,
and termed a {\em homoclinic corner} \cite{Si16b}.
An explicit expression for the line $\eta_3$ is available
because the left-most point of the chaotic attractor is $(-1,0)$,
while the left-most point of the $L^2 R$-cycle is $\left( x^{LRL}, 0 \right)$, where
\begin{equation}
x^{LRL} = \frac{\tau_L \tau_R + \tau_L - \delta_R + 1}{\tau_L^2 \tau_R - \tau_L \delta_R - 1}.
\nonumber
\end{equation}
The line $\eta_3$ is where these points coincide.
Solving $x^{LRL} = -1$ for $\delta_R$ gives
\begin{equation}
\delta_R = \tau_L \tau_R + \frac{\tau_L}{\tau_L + 1}.
\label{eq:eta3}
\end{equation}

In the grey region above $\eta_3$, the $L R^2$-cycle is only attractor.
At the left boundary $\alpha_3$ the $L R^2$-cycle tends to infinity,
while at the right boundary $\beta_3$ the $L R^2$-cycle loses stability
by attaining a stability multiplier of $-1$.
Formulas for these curves are derived in \S\ref{sub:neg-LRpm1}.

Upon crossing $\beta_3$ the attractor is again chaotic, see panel (a).
If we move further to the right the attractor destroyed at $\eta_4$
where it collides with an $L^2 R^2$-cycle.
Alternatively if we move down, then when we cross back over the line $\eta_3$
the attractor dramatically decreases in size, see panel (b).

\subsection{$L R^{p-1}$-regions}
\label{sub:neg-LRpm1}

Numerical explorations suggest that for all $p \ge 2$ there exists a region in Fig.~\ref{fig:zRbifSetNegative} where
\eqref{eq:f} has an asymptotically stable $L R^{p-1}$-cycle.
As $p$ increases these regions become narrower and exist for larger values of $\tau_R$.
Boundaries of the regions have been plotted for $p = 2,3,4,5$.

Here we derive formulas for these boundaries.
To do this we diagonalise $A_R$.
If $\delta_R \ne \frac{\tau_R^2}{4}$ then $A_R$ has distinct eigenvalues $\lambda_1$ and $\lambda_2$,
given by the roots of $\lambda^2 - \tau_R \lambda + \delta_R$.
We use these to form the matrices
\begin{equation}
D = \begin{bmatrix} \lambda_1 & 0 \\ 0 & \lambda_2 \end{bmatrix}, \qquad
P = \begin{bmatrix} 1 & 1 \\ -\lambda_2 & -\lambda_1 \end{bmatrix},
\nonumber
\end{equation}
giving the diagonalisation $A_R = P D P^{-1}$.
By evaluating $A_R^n = P D^n P^{-1}$ we obtain
\begin{equation}
A_R^n = \frac{1}{\lambda_1 - \lambda_2} \begin{bmatrix}
\lambda_1^{n+1} - \lambda_2^{n+1} & \lambda_1^n - \lambda_2^n \\
-\lambda_1 \lambda_2 \left( \lambda_1^n - \lambda_2^n \right) & \lambda_1 \lambda_2^n - \lambda_1^n \lambda_2
\end{bmatrix},
\label{eq:ARn}
\end{equation}
for any $n \in \mathbb{Z}$.
The stability multipliers of the $L R^{p-1}$-cycle are the eigenvalues of $M_{L R^{p-1}} = A_L A_R^{p-1}$.
This matrix has determinant zero and
\begin{equation}
{\rm trace} \left( M_{L R^{p-1}} \right) =
\frac{\tau_L \left( \lambda_1^p - \lambda_2^p \right) - \lambda_1 \lambda_2 \left( \lambda_1^{p-1} - \lambda_2^{p-1} \right)}
{\lambda_1 - \lambda_2},
\label{eq:traceALARpm1}
\end{equation}
obtained by using $n = p-1$ in \eqref{eq:ARn}.

The left boundaries of the $L R^{p-1}$ regions, labelled $\alpha_p$ for $p = 2,3$,
are where $M_{L R^{p-1}}$ has an eigenvalue $1$,
so are where ${\rm trace} \left( M_{L R^{p-1}} \right) = 1$.
In particular $\alpha_2$ is the line
\begin{equation}
\delta_R = \tau_L \tau_R - 1.
\label{eq:alpha2}
\end{equation}
The right boundaries, labelled $\beta_p$ for $p = 2,3$,
are where $M_{L R^{p-1}}$ has an eigenvalue $-1$,
so are where ${\rm trace} \left( M_{L R^{p-1}} \right) = -1$.
In particular $\beta_2$ is the line
\begin{equation}
\delta_R = \tau_L \tau_R + 1.
\label{eq:beta2}
\end{equation}
The bottom boundaries, labelled $\gamma_p$ for $p = 2,3$,
are BCBs where the point of the $L R^{p-1}$-cycle that maps into the left-half plane hits the switching manifold.
Straight-forward calculations lead to explicit formulas for these boundaries; for example $\gamma_2$ is given by
\begin{equation}
\delta_R = -\tau_R - 1,
\label{eq:gamma2}
\end{equation}
and $\gamma_3$ is given by
\begin{equation}
\delta_R = -\frac{(\tau_L + 1) \tau_R + 1}{\tau_L}.
\label{eq:gamma3}
\end{equation}

\subsection{Component doubling}

As we move about Fig.~\ref{fig:zRbifSetNegative}
the number of connected components of the chaotic attractor doubles or halves as we cross the dashed lines.
Specifically, as we move downwards across $\xi_1$ the
attractor changes from having one connected component, such as at the point {\bf b},
to having two connected components, such as at the point {\bf d}.
As we continue across $\xi_2$ the number of components increases to four,
then as we cross $\xi_3$ it increases to eight.
This phenomenon was identified by Kowalcyzk \cite[Proposition 2]{Ko05} up to four pieces.

To understand these lines, first observe that if $\delta_R = 0$ (the blue line in Fig.~\ref{fig:zRbifSetNegative})
then the dynamics reduces to one dimension.
Here every point maps to the $x$-axis
on which the dynamics is governed by the skew tent map \eqref{eq:skewTentMap}
with $(s_L,s_R,\eta) = (\tau_R,\tau_L,-\mu)$. 
So the blue line in Fig.~\ref{fig:zRbifSetNegative} (which uses $\tau_L = -1.2$)
corresponds to the horizontal line $s_R = -1.2$ in Fig.~\ref{fig:zRbifSetSkewTentMap}.

The dashed lines in Fig.~\ref{fig:zRbifSetNegative}
can be obtained by extending the construction of the dashed curves in Fig.~\ref{fig:zRbifSetSkewTentMap}
to accommodate $\delta_R \ne 0$.
For details of this type of construction refer to Ghosh {\em et al.}~\cite{GhSi22,GhMc24}.
The basic principle is that the attractor includes an area of phase space
where all points map under \eqref{eq:f} into the left-half plane $x<0$.
In this area the second iterate of the map is piecewise-linear with two pieces: $f_L^2$ and $f_L \circ f_R$.
The first piece has Jacobian matrix $A_L^2$,
which has determinant zero and trace $\tau_L^2$,
while the second piece has Jacobian matrix $A_L A_R$,
which has determinant zero and trace $\tau_L \tau_R - \delta_R$.
It follows that the $x$-dynamics is equivalent to the skew tent map
with $(s_L,s_R,\eta) = \left( \tau_L^2, \tau_L \tau_R - \delta_R, -\mu \right)$.
Thus we define
\begin{equation}
g_-(\tau_L,\tau_R,\delta_R) = \left( \tau_L^2, \tau_L \tau_R - \delta_R \right),
\label{eq:gMinus}
\end{equation}
and again use $\zeta$ given by \eqref{eq:divergenceFunction}.
The dashed line $\xi_1$ in Fig.~\ref{fig:zRbifSetNegative}
is $\zeta \left( g_-(\tau_L,\tau_R,\delta_R) \right) = 0$, so given by
\begin{equation}
\delta_R = \tau_L \tau_R + \frac{\tau_L^2}{\tau_L^2 - 1}.
\label{eq:zetaMinus1}
\end{equation}
The lines $\xi_2$ and $\xi_3$ are $\zeta \left( g^{k-1} \left( g_- (\tau_L,\tau_R,\delta_R) \right) \right) = 0$, for $k = 2$ and $k = 3$.
The next line in the sequence is not present for $\tau_L = -1.2$,
but using values of $\tau_L$ closer to $-1$
we can find lines where the number of components changes from $2^{k-1}$ to $2^k$ for all $k \ge 1$.

\subsection{Summary}

From the above observations we can make the following conclusions regarding
the long-term dynamics of \eqref{eq:f} with $\mu < 0$.

\begin{enumerate}
\renewcommand{\labelenumi}{\arabic{enumi})}
\item
Multiple attractors are possible, unlike in the one-dimensional setting.
In this case the long-term dynamics of typical orbits
depends critically on the initial point.
In the context of grazing-sliding bifurcations,
the creation of multiple attractors was reported in \cite{GlKo12,GlKo16,Si17d}.
\item
It appears that the only possible asymptotically stable periodic solutions are those with exactly one point in the left-half plane.
This is true for the one-dimensional map \eqref{eq:skewTentMap} with $\eta < 0$; a proof for the two-dimensional setting remains for future work.
The periodicity regions are ordered by period, hence a one-parameter bifurcation diagram corresponding to a path through
these regions will display {\em period-incrementing} \cite{AvGa19,DiBu08}.
\item
Chaotic attractors experience component doubling via the same mechanism as in one dimension.
Fig.~\ref{fig:zRbifSetNegative} uses $\tau_L = -1.2$ and chaotic attractors have $1$, $2$, $4$, or $8$ connected components.
Higher powers of two occur for values of $\tau_L$ closer to $-1$;
other numbers of components (e.g.~three) occur for large negative values of $\tau_L$ \cite{Ko05}.
\end{enumerate}

\section{Bifurcation structures for $\mu > 0$}
\label{sec:pos}

We now consider \eqref{eq:f} with $\mu > 0$.
Here the dynamics and bifurcations are more diverse
and we split this section into seven subsections corresponding to different aspects of the dynamics.

\begin{figure}[b!]
\begin{center}
\includegraphics[width=15.6cm]{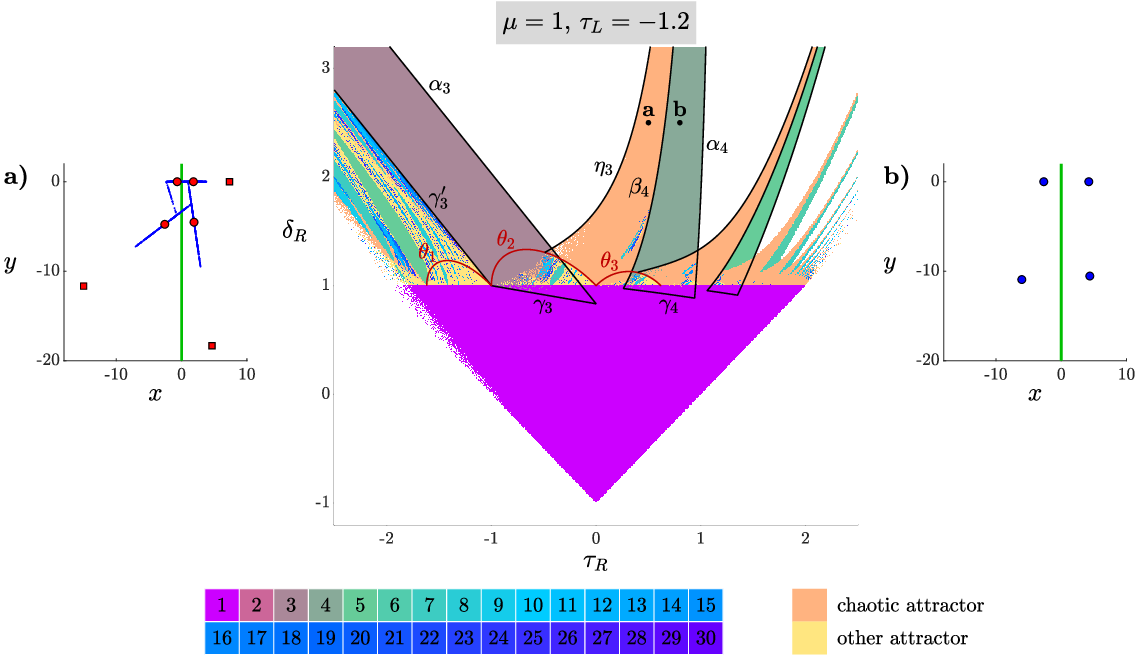}
\caption{
A two-parameter bifurcation diagram of \eqref{eq:f} with $\tau_L = -1.2$ and $\mu = 1$
computed by the numerical procedure of \S\ref{sub:numerics}.
Boundaries of regions where there exists an
asymptotically stable $L^2 R^{p-2}$-cycle are plotted for $p = 3,4,5$.
For $p = 3$ the boundaries are $\alpha_3$, $\gamma_3$, and $\gamma_3'$ given
by \eqref{eq:alpha3A}, \eqref{eq:gamma3A}, and \eqref{eq:gamma3pA} respectively.
The shrinking point curves $\theta_1$, $\theta_2$, and $\theta_3$ are given by \eqref{eq:theta1}, \eqref{eq:theta2}, and \eqref{eq:theta3} respectively.
The phase portraits correspond to the $(\tau_R,\delta_R)$ points
{\bf a}:~$(0.5,2.5)$;
{\bf b}:~$(0.8,2.5)$.
\label{fig:zRbifSetA}
} 
\end{center}
\end{figure}

To visualise parameter space we use two-parameter slices defined by fixing the value of $\tau_L$.
For any such slice the fixed point $(x^R,y^R)$ is admissible and asymptotically stable in 
the triangle $|\tau_R| - 1 < \delta_R < 1$.
Below we show slices at $\tau_L$-values of $-1.2$, $-0.4$, $0.4$, and $1.2$,
and the reader may first wish to glance ahead at these in
Figs.~\ref{fig:zRbifSetA}, \ref{fig:zRbifSetB}, \ref{fig:zRbifSetC}, and \ref{fig:zRbifSetD}.

\subsection{$L^2 R^{p-2}$-regions}

We start with Fig.~\ref{fig:zRbifSetA} which uses $\tau_L = -1.2$.
Here the dominant periodicity regions correspond to $L^2 R^{p-2}$ cycles.
The boundaries of these regions are shown in Figs.~\ref{fig:zRbifSetA} for $p = 3,4,5$.
Some of these boundaries are where the $L^2 R^{p-2}$ cycle has a stability multiplier of $1$ or $-1$.
The non-zero stability multiplier is the trace of $M_{L^2 R^{p-2}}$, given by
\begin{equation}
{\rm trace} \left( M_{L^2 R^{p-2}} \right) =
\frac{\tau_L^2 \left( \lambda_1^{p-1} - \lambda_2^{p-1} \right)
- \tau_L \lambda_1 \lambda_2 \left( \lambda_1^{p-2} - \lambda_2^{p-2} \right)}
{\lambda_1 - \lambda_2},
\label{eq:traceAL2ARpm2}
\end{equation}
using $n = p-2$ in \eqref{eq:ARn}.

For example at the parameter point {\bf b}, \eqref{eq:f} has an asymptotically stable $L^2 R^2$-cycle, see Fig.~\ref{fig:zRbifSetA}b.
As we move to the right the $L^2 R^2$-cycle tends to infinity when we cross $\alpha_4$
where ${\rm trace} \left( M_{L^2 R^2} \right) = 1$.
Alternatively as we move to left the $L^2 R^2$-cycle loses stability when we cross $\beta_4$
where ${\rm trace} \left( M_{L^2 R^2} \right) = -1$.
The third boundary of the $L^2 R^2$-region is $\gamma_4$ where one point of the $L^2 R^2$-cycle hits the switching manifold.

Immediately to the left of $\beta_4$, \eqref{eq:f} has a chaotic attractor.
Fig.~\ref{fig:zRbifSetA}a is a typical phase portrait showing the chaotic attractor (blue dots)
coexisting with the now unstable $L^2 R^2$-cycle (red circles)
as well as an unstable $L R^2$-cycle (red squares).
As we move further to the left the chaotic attractor is destroyed
when it collides with the $L R^2$-cycle on the curve $\eta_3$.

For periods greater than $p = 4$ the same four bifurcations occur at larger values of $\tau_R$.
In Fig.~\ref{fig:zRbifSetA} these bifurcation curves are plotted for $p = 5$.
With $p = 3$ there are analogous lines $\alpha_3$,
where ${\rm trace} \left( M_{L^2 R} \right) = 1$, and $\gamma_3$,
where one point of the $L^2 R$-cycle hits the switching manifold.
These lines are given by
\begin{equation}
\delta_R = \tau_L \tau_R - \frac{1}{\tau_L},
\label{eq:alpha3A}
\end{equation}
and
\begin{equation}
\delta_R = -\tau_R - \frac{1+\tau_R}{\tau_L},
\label{eq:gamma3A}
\end{equation}
respectively.
However, the left boundary $\gamma_3'$ of the $L^2 R$-region
is instead a BCB where a different point of the $L^2 R$-cycle hits the switching manifold;
this boundary is the line
\begin{equation}
\delta_R = \tau_L \tau_R + \tau_L + 1.
\label{eq:gamma3pA}
\end{equation}

\subsection{Centre bifurcations, sausage-strings, and period-adding}

The line $\delta_R = 1$ is a {\em centre bifurcation} where the fixed point $(x^R,y^R)$ loses stability.
This is a piecewise-linear analogue of a curve of Neimark-Sacker bifurcations
studied in detail by Sushko and Gardini \cite{SuGa08}
for the normal form \eqref{eq:bcnf} with $\delta_L \delta_R > 0$.
They found that immediately beyond the centre bifurcation there was often an attracting invariant circle
on which the dynamics was either quasi-periodic or mode-locked.
Mode-locked solutions occurred in Arnold tongues having a sausage-string geometry.

\begin{figure}[b!]
\begin{center}
\includegraphics[width=15.6cm]{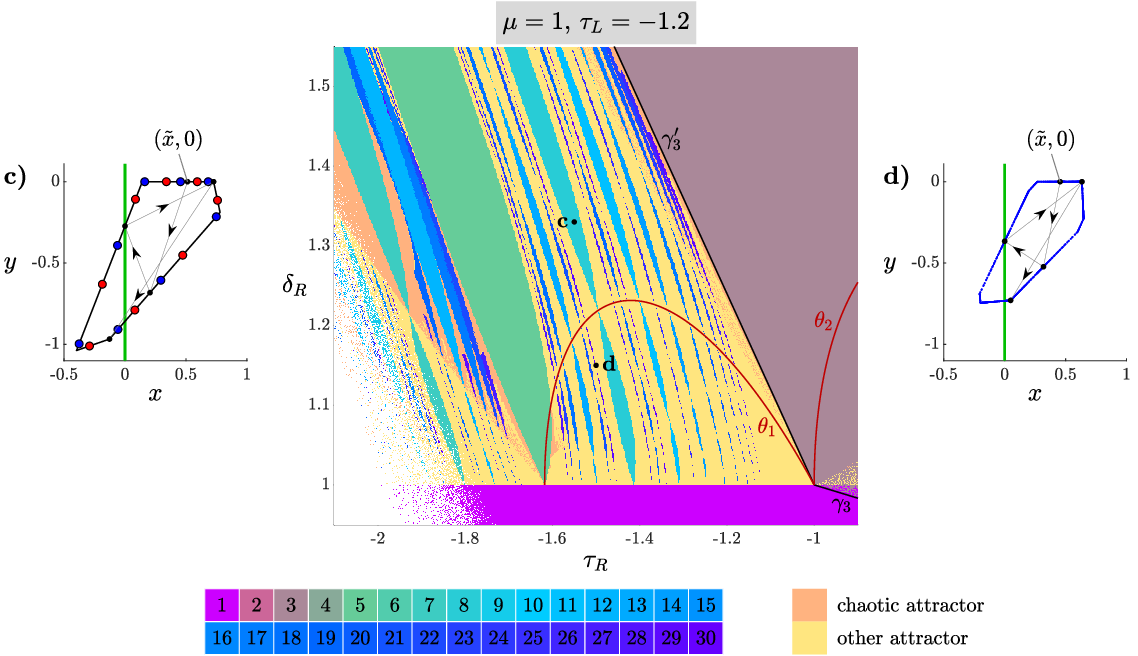}
\caption{
A magnification of Fig.~\ref{fig:zRbifSetA}.
The phase portraits correspond to the $(\tau_R,\delta_R)$ points
{\bf c}:~$(-1.55,1.33)$;
{\bf d}:~$(-1.5,1.15)$.
\label{fig:zRbifSetAZoom}
} 
\end{center}
\end{figure}

We observe the same behaviour here in the case $\delta_L = 0$,
at least along part of the $\delta_R = 1$ line.
This is shown more clearly for the magnification in Fig.~\ref{fig:zRbifSetAZoom}. 
For example the invariant circle is mode-locked to period eight at parameter point {\bf c}, shown in panel (c).
Here the circle contains an asymptotically stable $LRRLRRLR$-cycle (blue circles) and a saddle $RRRLRRLR$-cycle (red circles).
In contrast, at parameter point {\bf d} the dynamics on the invariant circle is either quasi-periodic
or mode-locked to a high period.

The average rate at which iterates step around the invariant circle is the {\em rotation number} \cite{AlLl00,DeVa93,KaHa95}.
Arnold tongues are regions where this number is constant and rational.
These regions emanate from $\delta_R = 1$ at $\tau_R$-values given by
\begin{equation}
\tau_R = 2 \cos(2 \pi \rho),
\label{eq:rotNum}
\end{equation}
where $\rho \in \mathbb{Q}$ is the rotation number.
This is because as we approach $\delta_R = 1$ the fraction of the invariant circle that belongs to the left-half plane tends to zero.
On $\delta_R = 1$ the fixed point $(x^R,y^R)$ has stability multipliers ${\rm e}^{\pm 2 \pi {\rm i} \rho}$, for some $\rho \in \mathbb{R}$,
so iterates in the right-half plane rotate around $(x^R,y^R)$ with rate $\rho$.
The stability multipliers sum to the trace $\tau_R$, hence obey \eqref{eq:rotNum}.
For example the blue Arnold tongue containing the point {\bf c}
has $\rho = \frac{3}{8}$, hence emanates from $\delta_R = 1$
at $\tau_R = 2 \cos \left( \frac{3 \pi}{4} \right) = -\sqrt{2}$.

As we move along the Arnold tongues by increasing the value of $\delta_R$,
some of the points in the corresponding stable periodic solutions move from the right-half plane to the left-half plane.
This cannot occur in a codimension-one fashion \cite{Si24d},
so the Arnold tongues have codimension-two points, termed {\em shrinking points},
where the width of the tongues shrinks to zero.
Consequently the Arnold tongues repeatedly widen then narrow giving them
a shape that is often likened to a string of sausages \cite{YaHa87}.
General theory for the bifurcations associated with the sausage-string structure can be found in \cite{Si17c,Si18e,SiMe09}.

The sausage-string structure was also studied by Szalai and Osinga \cite{SzOs09}
for the normal form \eqref{eq:bcnf} with $\delta_R = 0$.
The main difference to the general case of $\delta_L \ne 0$ and $\delta_R \ne 0$
is that the invariant circles always have the geometry of a polygon because they contain a segment of the $x$-axis
and can generated by taking finitely many images of this segment \cite{SzOs08}.


One-parameter bifurcation diagrams, say by taking a path through Fig.~\ref{fig:zRbifSetAZoom},
often cross several Arnold tongues resulting in intervals of periodicity.
As with smooth dynamical systems that have a persistent invariant circle,
lower periods typically correspond to wider periodicity intervals.
Furthermore, the rotation number varies continuously and often monotonically,
in which case the intervals display a {\em period-adding structure}.
That is, between intervals of periods $p_1$ and $p_2$,
the widest interval usually has period $p_1 + p_2$.
The main novelty of the piecewise-linear setting
is that the presence of shrinking points creates
a high degree of irregularity to the widths of the intervals.
This is seen in \S\ref{sec:appl} for the flu epidemic model.

\subsection{Shrinking point curves}

The red curve $\theta_1$ in Fig.~\ref{fig:zRbifSetAZoom} intersects Arnold tongues
at shrinking points and is an example of a {\em shrinking point curve}.
Such curves are in general where one point at which the invariant circle intersects the switching manifold
maps to the other point at which the invariant circle intersects the switching manifold in a fixed number of iterations.
For the normal form \eqref{eq:bcnf} with $\delta_L \ne 0$ and $\delta_R \ne 0$,
it was found in \cite{SiMe08b} that shrinking point curves are non-differentiable at shrinking points
because the position of the invariant circle varies with parameters in a complicated way.
But with $\delta_L = 0$ invariant circles contains an interval of the $x$-axis
and shrinking point curves are smooth and admit closed-form expressions.
This also occurs for families of piecewise-linear circle maps \cite{YaHa87}.

Let us now define $\theta_1$ and derive a formula for this curve.
The point $(\tilde{x},0)$, where
\begin{equation}
\tilde{x} = \frac{\tau_R + 1}{\delta_R - \tau_R^2},
\label{eq:xTilde}
\end{equation}
maps in two iterations to the switching manifold, see panels (a) and (b) of Fig.~\ref{fig:zRbifSetAZoom}.
If the map has an invariant circle containing $(\tilde{x},0)$,
then this circle contains the forward orbit of $(\tilde{x},0)$.
In this case the second iterate of $(\tilde{x},0)$ is one point at which the invariant circle intersects the switching manifold.
In panel (c) the fourth iterate of $(\tilde{x},0)$ lies to the left of the switching manifold,
while in panel (d) the fourth iterate of $(\tilde{x},0)$ lies to the right of the switching manifold.
We define $\theta_1$ to be where the fourth iterate of $(\tilde{x},0)$ lies on the switching manifold.
So here one intersection point of the invariant circle with the switching manifold maps to the other such point in two iterations.
By direct calculations we find $\theta_1$ is given implicitly by
\begin{equation}
\tau_R^3 + \tau_R^2 \delta_R + \tau_R \delta_R^2 + \tau_R^2 - \tau_R \delta_R - \delta_R = 0.
\label{eq:theta1}
\end{equation}
Notice \eqref{eq:theta1} is independent of $\tau_L$ because the four iterations only concern points with $x \ge 0$.

Now let $f$ denote the map \eqref{eq:f} and let us
characterise $\theta_1$ in a way that can be generalised to other shrinking point curves.
The curve $\theta_1$ is where there exists $x \in \mathbb{R}$ such that
$f^j(x,0)$ and $f^k(x,0)$ belong to the switching manifold, using $(j,k) = (2,4)$.
By using instead $(j,k) = (2,3)$ we obtain
\begin{equation}
\tau_R^2 + \tau_R \delta_R + \delta_R^2 - \delta_R = 0,
\label{eq:theta2}
\end{equation}
which is the curve $\theta_2$.
In Fig.~\ref{fig:zRbifSetA} it can be seen that the right part of this curve is a boundary of the chaotic region (orange).
Hence this is an example for which a shrinking point curve is boundary for chaos.
This phenomenon has been described for the normal form \eqref{eq:bcnf} with $\delta_L \ne 0$ and $\delta_R \ne 0$, \cite{SiMe08b}.

With instead $(j,k) = (3,4)$ we obtain
\begin{equation}
\tau_R^3 + \tau_R^2 \delta_R + \tau_R \delta_R^2 + \delta_R^3 - 2 \tau_R \delta_R - \delta_R^2 = 0,
\label{eq:theta3}
\end{equation}
which is the curve $\theta_3$.
In Fig.~\ref{fig:zRbifSetA} we see that the right part of this curve is also a boundary for chaos.
More such curves could be computed but we do not do this here.

\begin{figure}[b!]
\begin{center}
\includegraphics[width=15.6cm]{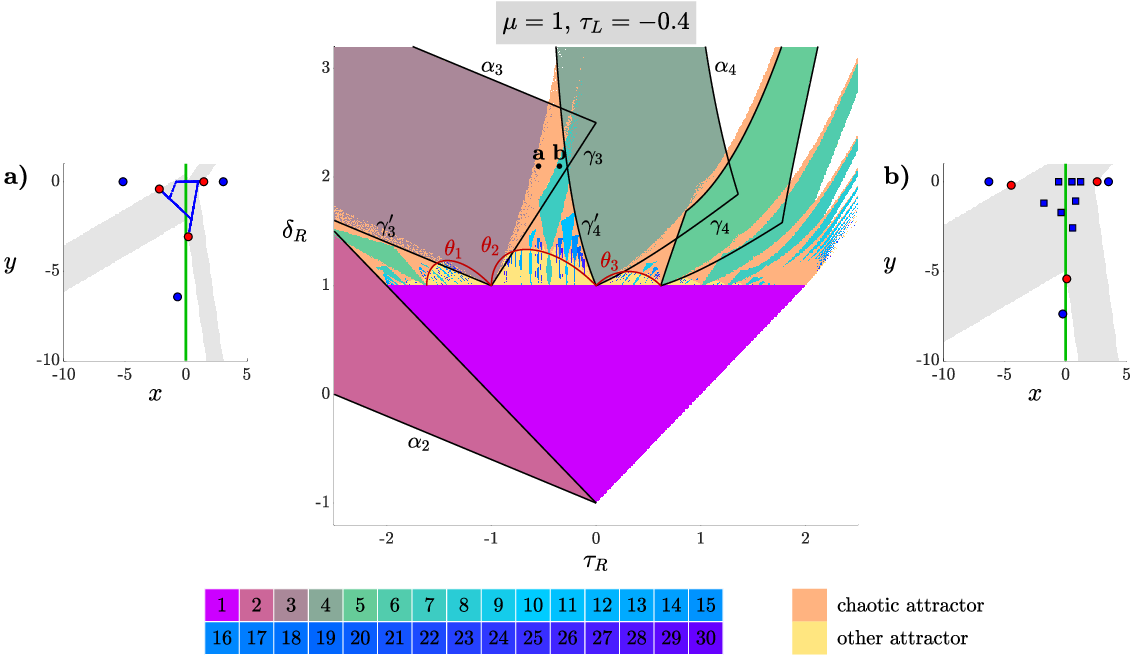}
\caption{
A two-parameter bifurcation diagram of \eqref{eq:f} with $\tau_L = -0.4$ and $\mu = 1$
computed by the numerical procedure of \S\ref{sub:numerics}.
We have plotted the boundaries of the region where there exists an asymptotically stable $LR$-cycle,
and the regions where there exists an
asymptotically stable $L^2 R^{p-2}$-cycle, for $p = 3,4,5$.
The boundaries $\alpha_p$ are where the periodic solution tends to infinity;
the boundaries $\gamma_p$ and $\gamma_p'$ are BCBs where one point of the periodic solution hits the switching manifold.
The phase portraits correspond to the $(\tau_R,\delta_R)$ points
{\bf a}:~$(-0.55,2.1)$, repeating Fig.~\ref{fig:zRbifDiagSchem}, and {\bf b}:~$(-0.35,2.1)$.
In (a) the basin of attraction of a chaotic attractor is shaded grey;
in (b) the basin of attraction of an $LLRRLRR$-cycle is shaded grey.
\label{fig:zRbifSetB}
} 
\end{center}
\end{figure}

Now consider Fig.~\ref{fig:zRbifSetB} which uses $\tau_L = -0.4$ instead of $\tau_L = -1.2$.
The location of the shrinking point curves $\theta_1$, $\theta_2$, and $\theta_3$
is unchanged, as they are independent of $\tau_L$, but now more of these curves are boundaries for chaos.

As before the dominant periodicity regions correspond to $L^2 R^{p-2}$-cycle.
Boundaries of these regions are shown for $p = 3,4,5$;
for $p = 3$ the boundaries are again
given by \eqref{eq:alpha3A}, \eqref{eq:gamma3A}, and \eqref{eq:gamma3pA}.
Note, in the lower-left part of Fig.~\ref{fig:zRbifSetB} there is now a dark pink region corresponding to an asymptotically stable $LR$-cycle.
The phase portraits in Fig.~\ref{fig:zRbifSetB} use parameter values at which the map has two attractors.
In both cases the stable manifold of a saddle $LR^2$-cycle bounds their basins of attraction.

\begin{figure}[b!]
\begin{center}
\includegraphics[width=15.6cm]{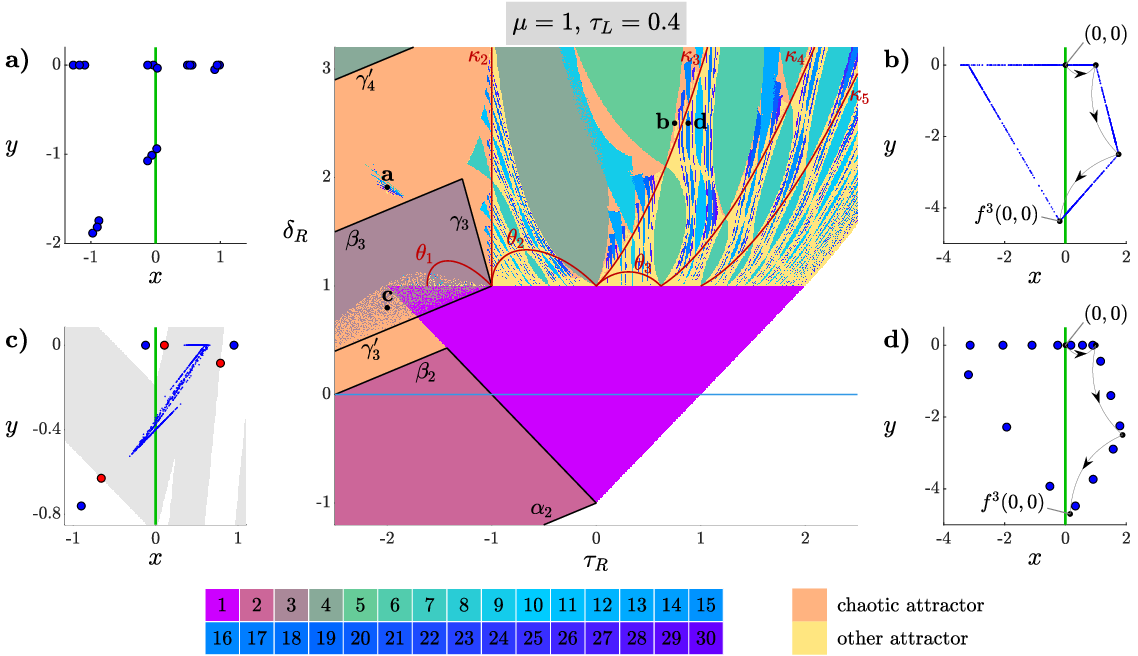}
\caption{
A two-parameter bifurcation diagram of \eqref{eq:f} with $\tau_L = 0.4$ and $\mu = 1$
computed by the numerical procedure of \S\ref{sub:numerics}.
Boundaries of regions where there exists an
asymptotically stable $L^{p-1} R$-cycle are plotted for $p = 2,3,4$.
The shrinking point curves $\theta_1$, $\theta_2$, and $\theta_3$
are given by \eqref{eq:theta1}, \eqref{eq:theta2}, and \eqref{eq:theta3},
while each $\kappa_n$ is given by \eqref{eq:shrPointCurves}.
The phase portraits correspond to the $(\tau_R,\delta_R)$ points
{\bf a}:~$(-2,1.91)$;
{\bf b}:~$(2.5,0.75)$;
{\bf c}:~$(-2,0.8)$;
{\bf d}:~$(2.5,0.88)$.
In (c) the basin of attraction of a chaotic attractor is shaded grey.
\label{fig:zRbifSetC}
} 
\end{center}
\end{figure}

Next Fig.~\ref{fig:zRbifSetC} uses $\tau_L = 0.4$.
Now the sausage-string geometry of the Arnold tongues extends higher, i.e.~to larger values of $\delta_R$,
and we have plotted additional shrinking point curves, $\kappa_2, \ldots, \kappa_5$.
Each $\kappa_n$ is defined to be where the origin maps to the switching manifold under $n$ iterations of the right piece of \eqref{eq:f}.
For instance at parameter point {\bf b} the third iterate of $(0,0)$ lies to the left of the switching manifold,
while at {\bf d} the third iterate of $(0,0)$ lies to the right of the switching manifold.
Thus $\kappa_3$ passes between {\bf b} and {\bf d}.
Note, {\bf d} belongs to an Arnold tongue with rotation number $\rho = \frac{3}{16}$,
so here the map has a stable period-$16$ solution,
while at {\bf b} the map appears to have a weakly chaotic attractor
including points that have `leaked' out of the invariant circle; compare \cite{SzOs08}.

Assuming \eqref{eq:f} has an invariant circle containing the origin,
$\kappa_n$ is where one point of intersection of the circle with the switching manifold maps to the other such point in $n$ iterations.
Notice the $\kappa_n$ intersect Arnold tongues exclusively at shrinking points.
Some of the $\kappa_n$ form a boundary of chaos; with $\tau_L = 0.4$ in Fig.~\ref{fig:zRbifSetC} this is particularly the case along $\kappa_2$.

Formulas for the $\kappa_n$ are straight-forward to derive.
For instance, $(0,0)$ maps to $(1,0)$ and then to $(\tau_R+1,-\delta_R)$ (using $\mu = 1$),
hence $\kappa_2$ is the line $\tau_R = -1$.
To handle larger values of $n$, observe
\begin{equation}
f_R^n(x,y;1) = A_R^n \begin{bmatrix} x \\ y \end{bmatrix} + \left( I + A_R + \cdots + A_R^{n-1} \right) \begin{bmatrix} 1 \\ 0 \end{bmatrix},
\nonumber
\end{equation}
and so
\begin{equation}
f_R^n(0,0;1) = \left( I - A_R \right)^{-1} \left( I - A_R^n \right) \begin{bmatrix} 1 \\ 0 \end{bmatrix}.
\nonumber
\end{equation}
By using \eqref{eq:ARn} we find that the first component of $f_R^n(0,0;1)$ is zero when
\begin{equation}
\lambda_1 - \lambda_2 - \left( \lambda_1^{n+1} - \lambda_2^{n+1} \right)
+ \lambda_1 \lambda_2 \left( \lambda_1^n - \lambda_2^n \right) = 0.
\label{eq:shrPointCurves}
\end{equation}
Thus $\kappa_n$ is given implicitly by \eqref{eq:shrPointCurves},
where $\lambda_1$ and $\lambda_2$ are the eigenvalues of $A_R$.

Additional shrinking point curves could be plotted, but we do not wish to clutter the figures.
Szalai and Osinga \cite{SzOs08} find that these curves are also where the number of sides
of the polygonal invariant circle changes.

\subsection{$L^{p-1} R$-regions}

\begin{figure}[b!]
\begin{center}
\includegraphics[width=15.6cm]{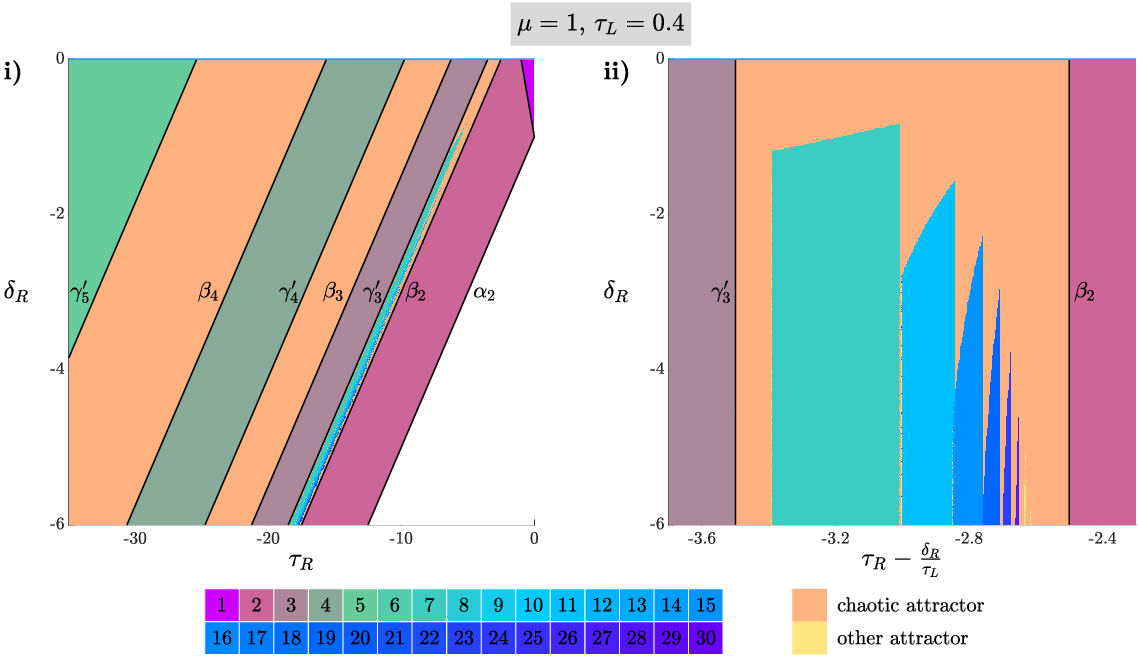}
\caption{
Panel (i) is the same as Fig.~\ref{fig:zRbifSetC}
but over a different range of $\tau_R$ and $\delta_R$ values.
Panel (ii) uses skewed coordinates to magnify the region between $\gamma_3'$ and $\beta_2$.
The tops of these plots are where $\delta_R = 0$
and the $x$-dynamics reduces to the skew tent map \eqref{eq:skewTentMap} with $(s_L,s_R,\eta) = (0.4,\tau_R,1)$.
\label{fig:zRbifSetCzoom}
} 
\end{center}
\end{figure}

Next we compute boundaries of parameter regions where there exist asymptotically stable $L^{p-1} R$-cycles.
These regions are present in the left part of Fig.~\ref{fig:zRbifSetC} for $p = 2,3,4$,
and seen more clearly in Fig.~\ref{fig:zRbifSetCzoom}(i) which includes large negative values of $\tau_R$,
and here part of the $p = 5$ region can also be seen.
The $L^{p-1} R$-regions are relatively large for small values of $\tau_L > 0$
(such as $\tau_L = 0.4$ in Figs.~\ref{fig:zRbifSetC} and \ref{fig:zRbifSetCzoom}).
This is because in the special case $\delta_R = 0$,
shown as blue lines in Figs.~\ref{fig:zRbifSetC} and \ref{fig:zRbifSetCzoom},
the $x$-dynamics of \eqref{eq:f}
is governed by the skew tent map \eqref{eq:skewTentMap} with $(s_L,s_R,\eta) = (\tau_L,\tau_R,\mu)$.
We saw in Fig.~\ref{fig:zRbifSetSkewTentMap}
that stable $L^{p-1} R$-cycles occur for the skew tent map with small $s_L > 0$, large $s_R < 0$, and any $\eta > 0$.

The $L^{p-1} R$-regions are ordered by their period $p$.
Thus one-parameter bifurcation diagrams that cross several of these regions will display period-incrementing
interspersed with bands of chaos.
As in the one-dimensional setting the lower boundaries $\beta_p$ of the $L^{p-1} R$-regions
are stability-loss boundaries where the $L^{p-1} R$-cycle has a stability multiplier of $-1$.
This stability multiplier is the trace of $A_R A_L^{p-1}$, which equals $\tau_L^{p-2} \left( \tau_L \tau_R - \delta_R \right)$.
By setting this equal to $-1$ and solving for $\delta_R$ we obtain
\begin{equation}
\delta_R = \tau_L \tau_R + \frac{1}{\tau_L^{p-2}},
\label{eq:pd5}
\end{equation}
for the lines $\beta_p$.

As in the one-dimensional setting the upper boundaries $\gamma_p'$ of the $L^{p-1} R$-regions
are BCBs where one point of the $L^{p-1} R$-cycle lies on the switching manifold.
To obtain formulas for these boundaries,
suppose a point $(x,0)$ maps under $f_R$ once, then under $f_L$ a total of $p-1$ times.
The resulting point is $(x',0)$, where
\begin{equation}
x' = \tau_L^{p-2} \left( \tau_L \tau_R - \delta_R \right) x + \frac{1 - \tau_L^p}{1 - \tau_L} \,\mu.
\label{eq:map5}
\end{equation}
The right-most point of the $L^{p-1} R$-cycle is $(x^*,0)$, where
\begin{equation}
x^* = \frac{1 - \tau_L^p}{(1 - \tau_L) \left( 1 - \tau_L^{p-2} (\tau_L \tau_R - \delta_R) \right)} \,\mu
\nonumber
\end{equation}
is the fixed point of \eqref{eq:map5}.
The BCB occurs when $x^* = \mu$ (so that its preimage is $(0,0)$ on the switching manifold),
and solving this equation for $\delta_R$ gives
\begin{equation}
\delta_R = \tau_L \tau_R - \frac{1}{\tau_L^{p-2}} \left( 1 - \frac{1 - \tau_L^p}{1 - \tau_L} \right).
\label{eq:bcb5}
\end{equation}
The lines \eqref{eq:bcb5} are labelled $\gamma_p'$ in Figs.~\ref{fig:zRbifSetC} and \ref{fig:zRbifSetCzoom}.

We finish this subsection with additional remarks regarding Figs.~\ref{fig:zRbifSetC} and \ref{fig:zRbifSetCzoom}.

\begin{remark}
As in the one-dimensional setting, neighbouring $L^{p-1} R$-regions are separated by bands of chaos,
but not between the $p=2$ and $p=3$ regions.
Fig.~\ref{fig:zRbifSetCzoom}(ii) uses skewed coordinates to blow up
the space between these regions which
contains a novel sequence of periodicity regions.
These correspond to stable $L^2 R (LR)^{2 k}$-cycles, for $k \ge 1$,
having periods $p = 3 + 4 k$.
An explanation for this structure remains for future work.
\end{remark}

\begin{remark}
Stable $L^{p-1} R$-cycles can coexist with other attractors.
For example at parameter point {\bf c} of Fig.~\ref{fig:zRbifSetC},
the stable $L^2 R$-cycle coexists with a chaotic attractor.
The stable manifold of a saddle $L R^2$-cycle (red circles in panel (c)) bounds their basins of attraction.
\end{remark}

\begin{figure}[b!]
\begin{center}
\includegraphics[width=12cm]{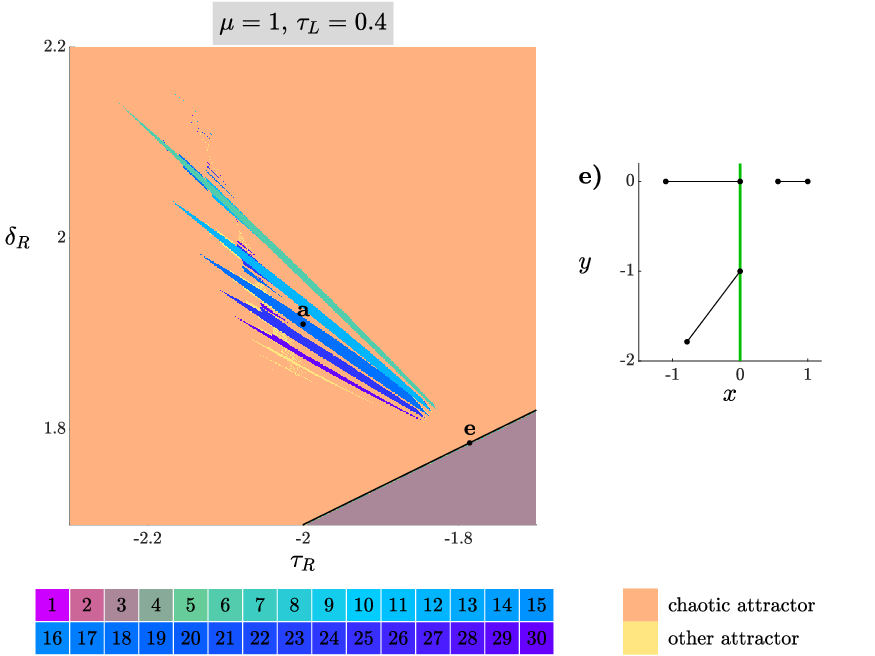}
\caption{
A two-parameter bifurcation diagram as in Fig.~\ref{fig:zRbifSetC}
except over a narrower range of $\tau_R$ and $\delta_R$ values.
The point {\bf a} is as in Fig.~\ref{fig:zRbifSetC};
the point {\bf e} is on $\beta_3$ with $\tau_R \approx -1.7857$
where the map has a period-$6$ solution with two points on the switching manifold.
The indicated line segments whose endpoints are the points of this periodic solution
consist entirely of $L^2 R L^2 R$-cycles.
\label{fig:zRbifSetCzoomAlt}
} 
\end{center}
\end{figure}

\begin{remark}
Above $\beta_3$ in Fig.~\ref{fig:zRbifSetC} the attractor is chaotic except
in a small group of regions roughly centred at parameter point {\bf a}
(here \eqref{eq:f} has a stable period-$18$ solution).
This is seen more clearly in the magnification provided by Fig.~\ref{fig:zRbifSetCzoomAlt}.
The periodicity regions emanate from the codimension-two point {\bf e}, which is on $\beta_3$
and where \eqref{eq:f} has an $L^2 R L^2 R$-cycle with two points on the switching manifold,
although there is a gap between the periodicity regions and {\bf e}.
This type of structure was reported in \cite[Figure 20]{SiMe08b},
but a theory that explains the structure has not been described.
\end{remark}

\subsection{Homoclinic corners}

\begin{figure}[b!]
\begin{center}
\includegraphics[width=15.6cm]{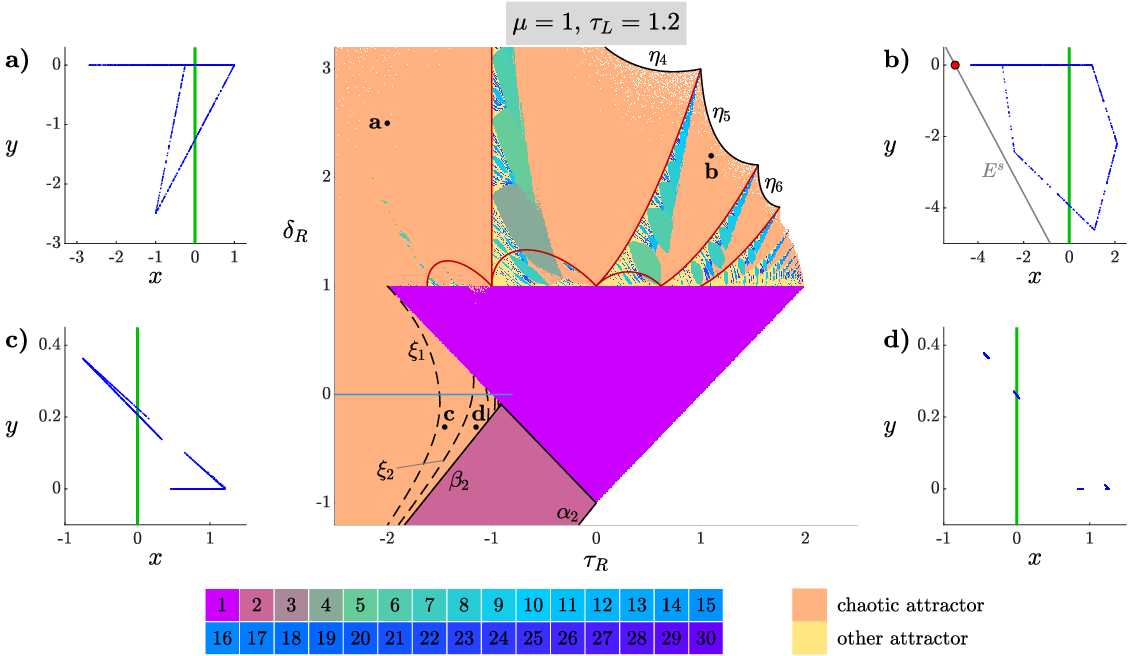}
\caption{
A two-parameter bifurcation diagram of \eqref{eq:f} with $\tau_L = 1.2$ and $\mu = 1$
computed by the numerical procedure of \S\ref{sub:numerics}.
We again show the shrinking point curves of Fig.~\ref{fig:zRbifSetC},
but now also the homoclinic corners $\eta_n$, for $n = 4,5,6$, discussed in the text.
The horizontal blue line at $\delta_R = 0$ corresponds
to the vertical blue line at $s_L = 1.2$ in the lower-right plot of Fig.~\ref{fig:zRbifSetSkewTentMap}.
The dashed curves $\xi_k$ are where the number of connected components of
the chaotic attractor doubles (or halves).
The phase portraits correspond to the $(\tau_R,\delta_R)$ points
{\bf a}:~$(-2,2.5)$;
{\bf b}:~$(1.1,2.2)$;
{\bf c}:~$(-1.45,-0.3)$;
{\bf d}:~$(-1.15,-0.3)$,
and show the attractor in blue.
In panel (b) we also show the saddle fixed point $(x^L,y^L)$
and the initial linear part of its stable manifold.
\label{fig:zRbifSetD}
} 
\end{center}
\end{figure}

The last value of $\tau_L$ that we examine is $\tau_L = 1.2$, Fig.~\ref{fig:zRbifSetD}.
In comparison to earlier figures the only remaining $L^{p-1} R$-region has $p=2$.
There are still Arnold tongues,
but many of these are not connected to the centre bifurcation $\delta_R = 1$.
Chaos is more prevalent; almost the entire visible length of the shrinking point curves $\kappa_2, \ldots, \kappa_5$
are boundaries between chaotic and non-chaotic attractors.

The boundaries from orange to white in the top-right of Fig.~\ref{fig:zRbifSetD}
are where a chaotic attractor is destroyed, analogous to the line $\eta_3$ of Fig.~\ref{fig:zRbifSetNegative} in the case $\mu < 0$.
Boundaries similar to those in Fig.~\ref{fig:zRbifSetD}
have been reported for the normal form \eqref{eq:bcnf} with $\delta_L \ne 0$ and $\delta_R \ne 0$, see \cite{GlSi22b,Si16b,Si20}.
In general these are homoclinic corners where the stable and unstable manifolds of $(x^L,y^L)$ attain a non-trivial intersection.
This is because the attractor is contained in the closure of the right branch of the unstable manifold of $(x^L,y^L)$,
and as soon as this manifold intersects the stable manifold of $(x^L,y^L)$ at points other than $(x^L,y^L)$,
there is a route for orbits on the right branch of the unstable manifold to diverge.
Note that $(x^L,y^L)$ is a saddle throughout Fig.~\ref{fig:zRbifSetD} because $\tau_L > 1$.

The novelty here is that $\delta_L = 0$, so the homoclinic corners can be characterised
by simply where $(0,0)$ maps to $(x^L,y^L)$ in a few iterations.
To explain, Fig.~\ref{fig:zRbifSetD}b shows the chaotic attractor at parameter point {\bf b}
located just below the homoclinic corner $\eta_5$.
In this phase portrait the fixed point $(x^L,y^L)$ is indicated with a red circle,
while the initial linear part of its stable manifold $E^s$ is shown in grey.
The unstable manifold of $(x^L,y^L)$ emanates from $(x^L,y^L)$ along the $x$-axis
and has corners at the images of $(0,0)$.
The chaotic attractor is contained in the unstable manifold of $(x^L,y^L)$,
and its left-most point is the fifth iterate of $(0,0)$.
Each $\eta_n$ in Fig.~\ref{fig:zRbifSetD}
is where the $n^{\rm th}$ iterate of $(0,0)$ under \eqref{eq:f} equals $(x^L,y^L)$.

\subsection{Superstable periodic solutions}

\begin{figure}[b!]
\begin{center}
\includegraphics[width=15.6cm]{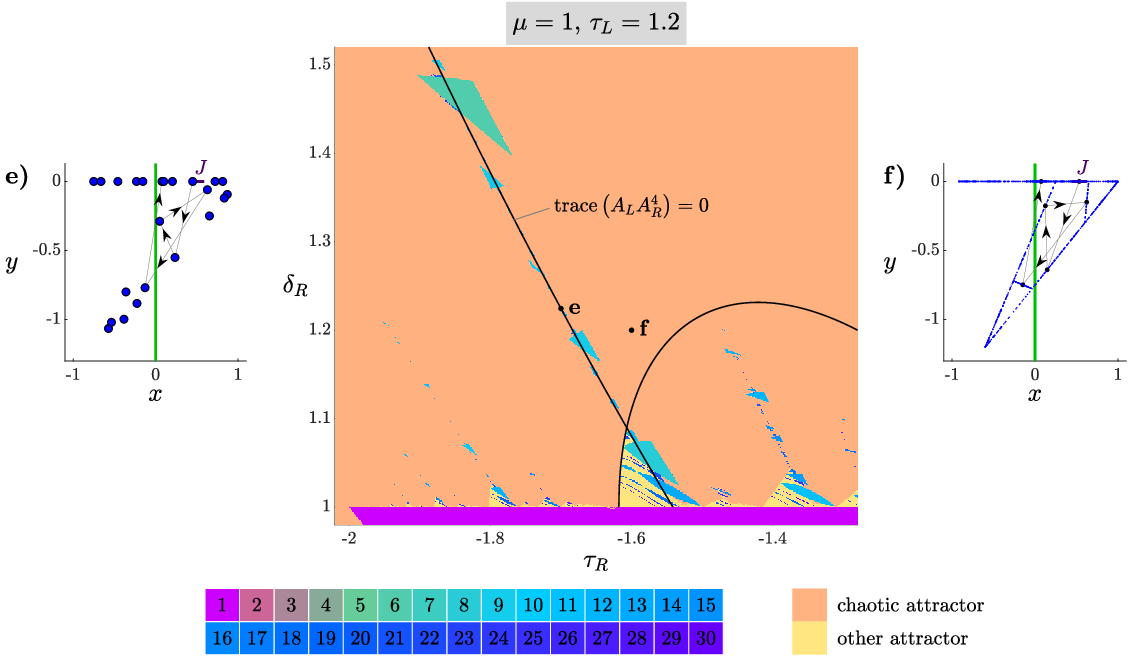}
\caption{
A two-parameter bifurcation diagram using the same parameter values as Fig.~\ref{fig:zRbifSetD}
but over a smaller range of $\tau_R$ and $\delta_R$ values.
The phase portraits correspond to the $(\tau_R,\delta_R)$ points
{\bf e}:~$(-1.7,\delta_R^*)$, where $\delta_R^* \approx 1.2245$
is such that ${\rm trace} \left( A_L A_R^4 \right) = 0$, and
{\bf f}:~$(-1.6,1.2)$.
The interval $J$ has endpoints at $x = -\frac{1}{\tau_R}$ and $x = \tilde{x}$, given by \eqref{eq:xTilde},
and consists of all points $(x,0)$
for which $f^i(x,0)$ belongs to the right-half plane for all $i = 0,1,2,3$ and to the left-half plane for $i = 4$,
where $f$ is the map \eqref{eq:f}.
\label{fig:zRbifSetDzoom}
} 
\end{center}
\end{figure}

Fig.~\ref{fig:zRbifSetDzoom} is a magnification of Fig.~\ref{fig:zRbifSetD}
and shows several isolated periodic regions that do not display the sausage-string structure.
These regions are mostly triangular in shape and occur roughly along curves emanating from $\delta_R = 1$.
This type of structure does not seem to have been reported before.

We have noticed that the largest sequence of these regions
occur along the curve defined by ${\rm trace} \left( A_L A_R^4 \right) = 0$ and given by
\begin{equation}
\left( \tau_L + 2 \tau_R \right) \delta_R^2 - \left( 3 \tau_L + \tau_R \right) \tau_R^2 \delta_R
+ \tau_L \tau_R^4 = 0.
\label{eq:mCzero}
\end{equation}
To explain the significance of this curve,
observe that for parameter combinations throughout most of Fig.~\ref{fig:zRbifSetDzoom}
there exists an interval $J$ on the $x$-axis
whose points map under the right piece of \eqref{eq:f} four times, then the left piece of \eqref{eq:f} once,
upon which they return to the $x$-axis, see panels (e) and (f) of Fig.~\ref{fig:zRbifSetDzoom}.
But if ${\rm trace} \left( A_L A_R^4 \right) = 0$ then all points in $J$ return to the same point on the $x$-axis.
In this case any periodic solution of \eqref{eq:f} having a point in $J$ is superstable:
both of its stability multipliers are zero.
In this way the curve ${\rm trace} \left( A_L A_R^4 \right) = 0$ inhibits the occurrence of a chaotic attractor.

It remains to characterise other curves emanating from $\delta_R = 1$, and determine where and how these terminate.
It would also be helpful to explain what governs the size, shape, and ordering of the periodicity regions along these curves.

\subsection{Component doubling}

With $\delta_R = 0$ in Fig.~\ref{fig:zRbifSetD}
the $x$-dynamics of \eqref{eq:f} are governed by the skew tent map \eqref{eq:skewTentMap}
with $(s_L,s_R,\eta) = (1.2,\tau_R,1)$.
As seen in the lower-right plot of Fig.~\ref{fig:zRbifSetSkewTentMap},
these values include chaotic attractors with component doubling as the value of $s_R$ is varied.

Fig.~\ref{fig:zRbifSetD} shows how the component doubling boundaries extend to $\delta_R \ne 0$
and were obtained through renormalisation, similar to those in Fig.~\ref{fig:zRbifSetNegative}.
Panels (a), (b), and (d)
correspond to points on different sides of these curves
and show chaotic attractors having one, two, and four connected components, respectively.

\subsection{Summary}

From the above observations we can make the following conclusions regarding
the long-term dynamics of \eqref{eq:f} with $\mu > 0$.

\begin{enumerate}
\renewcommand{\labelenumi}{\arabic{enumi})}
\item
As with $\mu < 0$, multiple attractors are possible.
\item
As with $\mu < 0$, period-incrementing structures occur, see Fig.~\ref{fig:zRbifSetCzoom}a.
But now period-adding structures also occur.
In two-parameter bifurcation diagrams periodicity regions have a sausage-string structure,
hence in one-parameter bifurcation diagrams the widths of periodicity intervals
can usually be expected to have a high degree of irregularity.
\item
Unlike for $\mu < 0$, quasi-periodic motion is now possible.
\item
Chaotic attractors occur in relatively large open regions of parameter space
and experience component doubling in a similar fashion to the one-dimensional setting.
\item
Chaotic attractors are often destroyed in homoclinic corners beyond which there is no attractor because typical forward orbits diverge.
They are also destroyed along shrinking point curves beyond which the attractor is periodic or quasi-periodic.
\end{enumerate}

\section{Applications}
\label{sec:appl}

In this section we illustrate the above results with three applications.
We revisit the stick-slip friction oscillator models studied by Yoshitake and Sueoka \cite{YoSu00}
and Szalai and Osinga \cite{SzOs08},
then study the epidemiology model of Roberts {\em et al.}~\cite{RoHi19b}.

\subsection{A nonlinear stick-slip friction oscillator}

\begin{figure}[b!]
\begin{center}
\includegraphics[width=9cm]{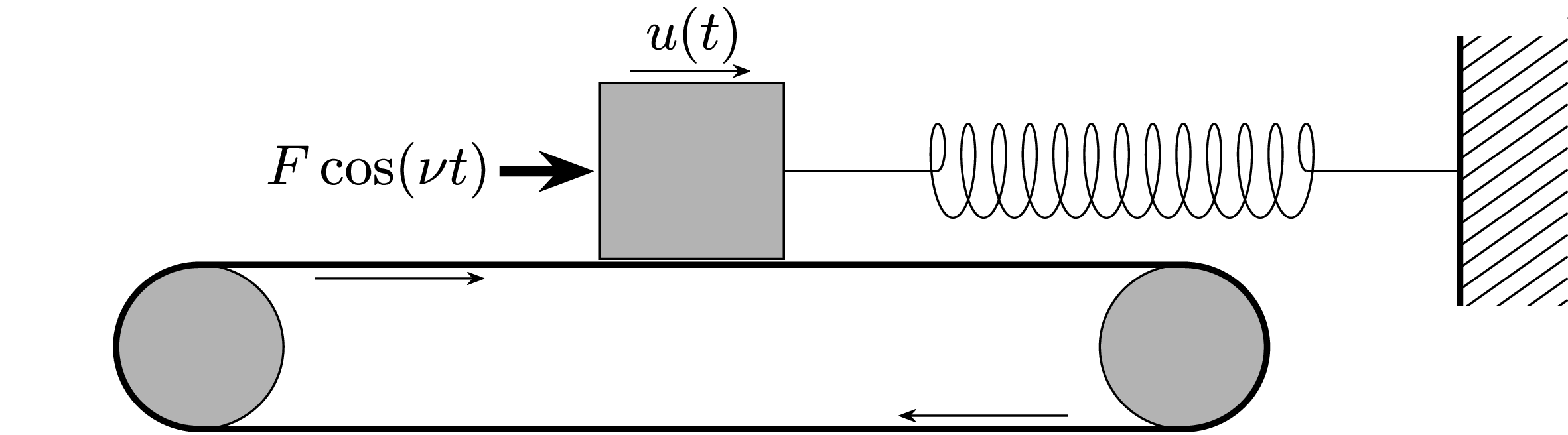}
\caption{
The block-belt system modelled by \eqref{eq:DiKo03ODEs}.
The horizontal displacement $u(t)$ of the block
is governed by a spring fixed to a wall,
a belt rotating at constant velocity,
and a external force of $F \cos(\nu t)$.
\label{fig:zRFrictionOscSchem}
} 
\end{center}
\end{figure}

For the block-belt system of Fig.~\ref{fig:zRFrictionOscSchem} we assume the position $u(t)$ of the block obeys
\begin{equation}
\ddot{u} + u = \alpha_0 \,{\rm sgn} \left( 1 - \dot{u} \right)
- \alpha_1 \left( 1 - \dot{u} \right)
+ \alpha_2 \left( 1 - \dot{u} \right)^3 + F \cos(\nu t),
\label{eq:DiKo03ODEs}
\end{equation}
where $\dot{u} = \frac{d u}{d t}$ and $\ddot{u} = \frac{d^2 u}{d t^2}$.
The belt rotates with constant velocity.
The block is at some times stuck to the belt, due to static friction with maximal force $\alpha_0 > 0$,
and at other times slips over the belt, whereby the magnitude of friction
is specified through the kinetic friction parameters $\alpha_1$ and $\alpha_2$.
With $\alpha_1 > 0$ and $\alpha_2 > 0$ the kinetic friction has negative slope for small $|1 - \dot{u}|$
and positive slope of large $|1 - \dot{u}|$.
This replicates the Stribeck effect whereby an extra force is required
for the block to transition from sticking to slipping \cite{LuKh06}.
The block is also attached to a spring and forced harmonically with amplitude $F$ and frequency $\nu$.

Stick-slip motion is realised by treating \eqref{eq:DiKo03ODEs} as a Filippov system \cite{DiBu08,Fi88}.
That is, $\dot{u} = 1$ is a discontinuity surface in phase space
on which sliding motion corresponds to the block being stuck to the belt.

\begin{figure}[b!]
\begin{center}
\includegraphics[width=15.6cm]{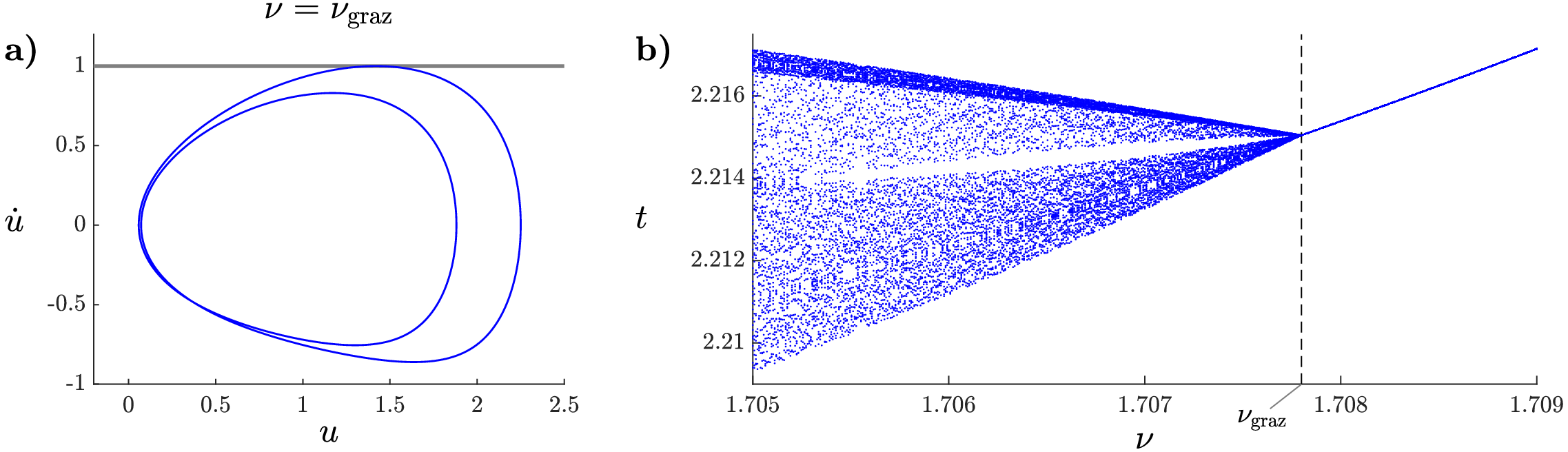}
\caption{
Panel (a) shows the stable limit cycle of \eqref{eq:DiKo03ODEs} with \eqref{eq:DiKo03params},
$F = 0.1$ and $\nu = \nu_{\rm graz} \approx 1.7078$.
Panel (b) is a bifurcation diagram showing how the attractor of \eqref{eq:DiKo03ODEs} changes with the value of $\nu$.
\label{fig:zRDiKo03_A}
} 
\end{center}
\end{figure}

As in di Bernardo {\em et al.}~\cite{DiKo03}, we fix
\begin{equation}
\alpha_0 = 1.5, \qquad
\alpha_1 = 1.5, \qquad
\alpha_2 = 0.45.
\label{eq:DiKo03params}
\end{equation}
and consider various values of $F$ and $\nu$.
With $F = 0.1$ the system undergoes a grazing-sliding bifurcation at $\nu = \nu_{\rm graz} \approx 1.7078$
by having an asymptotically stable limit cycle hit the discontinuity surface.
Fig.~\ref{fig:zRDiKo03_A}a shows the limit cycle at grazing (it has two loops and period $T = \frac{8 \pi}{\nu_{\rm graz}}$).
Fig.~\ref{fig:zRDiKo03_A}b is a bifurcation diagram showing that the attracting solution becomes chaotic
as the value of $\nu$ is decreased through $\nu_{\rm graz}$.
Initially the chaos is relatively weak as solutions oscillate irregularly but
close to the orbit shown in Fig.~\ref{fig:zRDiKo03_A}a.
As $\nu$ is decreased further the amplitude of the chaos increases steadily.

The vertical axis of Fig.~\ref{fig:zRDiKo03_A}b
shows times, modulo $\frac{2 \pi}{\nu}$, at which $\ddot{u} = 0$
for numerically computed solutions to \eqref{eq:DiKo03ODEs} with transient behaviour removed.
Notice that for values of $\nu$ just less than $\nu_{\rm graz}$
there is a gap amongst the scattering of times.
This indicates that, in the context of the Poincar\'e map with domain $\ddot{u} = 0$,
the chaotic attractor created at the grazing-sliding bifurcation has two connected components.

To investigate this more carefully,
we evaluate derivatives of the Poincar\'e map.
The Poincar\'e map is piecewise-smooth, with two $C^1$ pieces.
We used finite differences to approximate the derivative of each piece of the map at the bifurcation point for $\nu = \nu_{\rm graz}$.
We found trace and determinant of one derivative was $-1.653$ and $0$,
while the trace and determinant of the other derivative was $0.848$ and $0.006$, approximately.
Thus the grazing-sliding bifurcation corresponds to
the map \eqref{eq:f} with $(\tau_L,\tau_R,\delta_R) \approx (-1.653,0.848,0.006)$.
More precisely, the piecewise-linear approximation to the Poincar\'e map
is affinely conjugate to \eqref{eq:f} with these parameter values.

This instance of \eqref{eq:f} has a stable fixed point for $\mu > 0$
because $(\tau_R,\delta_R) \approx (0.848,0.006)$ satisfies $|\tau_R| - 1 < \delta_R < 1$
(the purple triangles in the earlier figures). 
To get some idea of the dynamics of \eqref{eq:f} for $\mu < 0$
we can refer to Fig.~\ref{fig:zRbifSetNegative} which uses $\tau_L = -1.2$.
With $\tau_L \approx -1.653$
the point $(\tau_R,\delta_R) \approx (0.848,0.006)$
lies between $\xi_1$ and $\xi_2$
so \eqref{eq:f} has a chaotic attractor with two connected components.
In this way the above theory explains the features of Fig.~\ref{fig:zRDiKo03_A}b.

\begin{figure}[b!]
\begin{center}
\includegraphics[width=15.6cm]{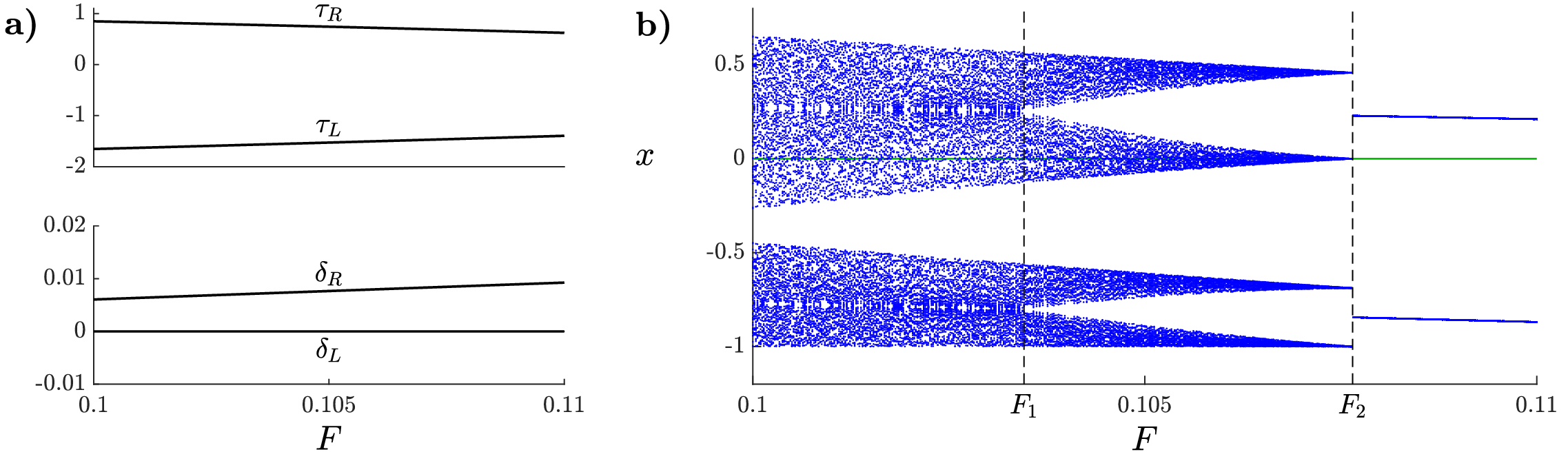}
\caption{
Panel (a) shows the parameter values of the normal form \eqref{eq:bcnf}
for the grazing-sliding bifurcation of \eqref{eq:DiKo03ODEs} with \eqref{eq:DiKo03params}.
Panel (b) indicates the attractor of \eqref{eq:f} with $\mu = 1$.
\label{fig:zRDiKo03_B}
} 
\end{center}
\end{figure}

For different values of the forcing amplitude $F$,
the grazing-sliding bifurcation occurs at different values of $\nu$
and corresponds to different values of $\tau_L$, $\tau_R$, and $\delta_R$.
This is indicated in Fig.~\ref{fig:zRDiKo03_B}a.
Fig.~\ref{fig:zRDiKo03_B}b shows how the attractor of the corresponding instance of \eqref{eq:f}
changes over the same range of values of $F$.
This represents the attractor created at the grazing-sliding bifurcation.

As $F$ increases, $\tau_L$ increases while $\tau_R$ decreases.
At $F = F_1 \approx 0.10346$ the point $(\tau_R,\delta_R)$ crosses $\xi_2$
beyond which the attractor has four connected components.
Then at $F = F_2 \approx 0.10765$ the point $(\tau_R,\delta_R)$ crosses $\beta_2$, given by \eqref{eq:beta2},
and beyond this the map has a stable $LR$-cycle.
In terms of the block-belt system
this corresponds to a period $\frac{16 \pi}{\nu}$ solution
with one short phase of sticking motion per period.

\subsection{A linear stick-slip friction oscillator}

Szalai and Osinga \cite{SzOs08} study a model equivalent to \eqref{eq:DiKo03ODEs} with $\alpha_2 = 0$.
Again grazing-sliding bifurcations occur, although now the grazing limit cycle is unstable.
Since there is now no cubic term,
the flow of \eqref{eq:DiKo03ODEs} can be expressed explicitly when $1 - \dot{u}$ has constant sign.
It follows that the values of $\tau_L$, $\delta_L$, $\tau_R$, and $\delta_R$
associated with the grazing-sliding bifurcation
can be expressed explicitly in terms of the model parameters.
Specifically,
\begin{equation}
\tau_L = \re^\beta \cos(\alpha), \quad
\delta_L = 0, \quad
\tau_R = 2 \re^\beta \cos(\alpha), \quad
\delta_R = \re^{2 \beta},
\label{eq:of}
\end{equation}
where $\alpha = \frac{2 \pi}{\nu} \sqrt{1 - \frac{\alpha_1^2}{4}}$
and $\beta = \frac{\pi \alpha_1}{\nu}$.
By applying these formulas to the results of the previous sections
we can explain the dynamics local to the grazing-sliding bifurcation.

\begin{figure}[b!]
\begin{center}
\includegraphics[width=15.6cm]{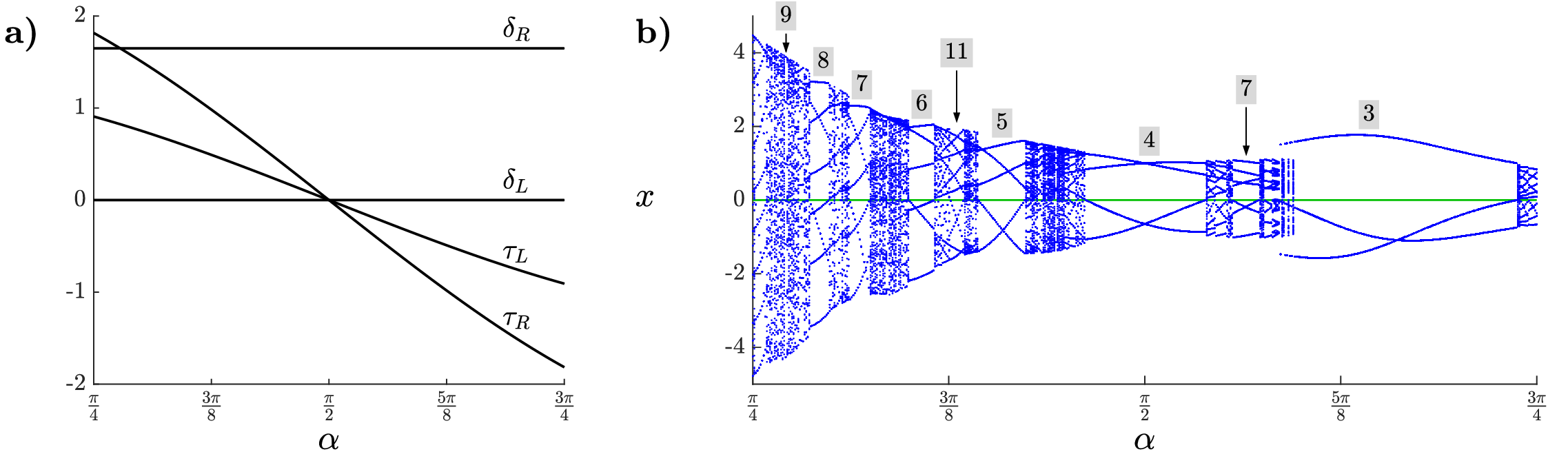}
\caption{
Panel (a) shows the values \eqref{eq:of} with $\beta = 0.25$.
Panel (b) indicates the attractor of \eqref{eq:f} with $\mu = 1$.
\label{fig:zRSzOs08_B}
} 
\end{center}
\end{figure}

So for example if we fix $\beta = 0.25$ and vary $\alpha$
we obtain the values of $\tau_L$, $\delta_L$, $\tau_R$, and $\delta_R$ shown in Fig.~\ref{fig:zRSzOs08_B}a.
Notice $|\tau_L| < 1$, for all values of $\alpha$ shown,
thus for $\mu = -1$ the fixed point $(x^L,y^L)$ of \eqref{eq:f} is admissible and asymptotically stable.
This corresponds to a stable limit cycle of the block-belt system
having one brief phase of sticking motion per period.

Fig.~\ref{fig:zRSzOs08_B}b provides a corresponding bifurcation diagram of \eqref{eq:f} with $\mu = 1$.
This shows how the attractor of the block-belt system
created in the grazing-sliding bifurcation varies with $\alpha$.
It shows period-adding, which is consistent with our earlier figures.
Here $\delta_R = \sqrt{e} \approx 1.6487$, while the value of $\tau_R$ is twice that of $\tau_L$.
Figs.~\ref{fig:zRbifSetB} and \ref{fig:zRbifSetC}
suggest this puts us in a regime dominated by periodicity regions with the sausage-string structure.
Fig.~\ref{fig:zRSzOs08_B}b represents a one-parameter path through these regions
and indeed displays periodic intervals with period-adding.
For example between the period-$3$ and period-$4$ intervals there is a period-$7$ interval,
and between the period-$5$ and period-$6$ intervals there is a period-$11$ interval.
As a consequence of the sausage-string structure the widths of the periodicity intervals is somewhat irregular.
For example the period-$9$ interval in the left of the figure is almost too narrow to be seen.
This is because the one-parameter path defined by fixing $\beta = 0.25$
passes close a shrinking point of an Arnold tongue with period-$9$,
see \cite[Figure 1]{Si20d}.

\subsection{Patterns of influenza outbreaks}

Here we study the influenza outbreak model of Roberts {\em et al.}~\cite{RoHi19b}.
This combines a standard SIR model for the evolution of influenza cases during flu seasons,
with a map that captures the change in the status of population from the end of one flu season
to the start of the next flu season.
The model assumes that during each flu season
the distribution of the population across the SIR compartments
reaches equilibrium by the end of the flu season.
This equilibrium has two different functional forms depending on whether or not an outbreak occurs.
Consequently the formula we extract from the model for
the status of the population at the start of a flu season
in terms of the status of the population at the start of the previous flu season
is a piecewise-smooth map.
This map is
\begin{equation}
\begin{bmatrix}
S \\ T
\end{bmatrix}
\mapsto \begin{bmatrix}
1 + c(S + T - 1 - p(S,T)) \\
-c(S + T - 1)
\end{bmatrix},
\label{eq:m10-RoHi19b}
\end{equation}
where $p = p(S,T)$ satisfies
\begin{equation}
p = \begin{cases}
0, & r(S,T) \le 1, \\
S \left( 1 - \re^{-R_0 p} \right) + T \left( 1 - \re^{-k R_0 p} \right), & r(S,T) \ge 1,
\end{cases}
\label{eq:m10-RoHi19b-p}
\end{equation}
and $r(S,T) = R_0 (S + k T)$.
The variables $S$ and $T$
represent the fraction of the population that is fully susceptible and partially susceptible to influenza, respectively.
The parameter $k \in (0,1)$ is the relative susceptibility of partially susceptible individuals compared to fully susceptible individuals,
$c \in (0,1)$ is the probability that an immune individual is still immune at the start of the next season,
and $R_0 > 0$ is the basic reproduction number.

If $r(S,T) < 1$ then no outbreak occurs.
In this case the image $(S',T')$ of \eqref{eq:m10-RoHi19b}
satisfies $S' + T' = 1$, meaning that at the start of next flu season the entire population is susceptible, i.e.~no-one is infected.
This also means that the $r(S,T) < 1$ piece of the map has degenerate range.
For this reason will be able to use \eqref{eq:f}, as explained below.

If $r(S,T) > 1$ then an outbreak does occur.
In this case the value of $p$, which represents the fraction of the population that is infected during the outbreak,
is given {\em implicitly} by the second piece of \eqref{eq:m10-RoHi19b-p}.
To implement \eqref{eq:m10-RoHi19b} numerically, we used a root-finding method to solve \eqref{eq:m10-RoHi19b-p} for $p$.

\begin{figure}[b!]
\begin{center}
\includegraphics[width=15.6cm]{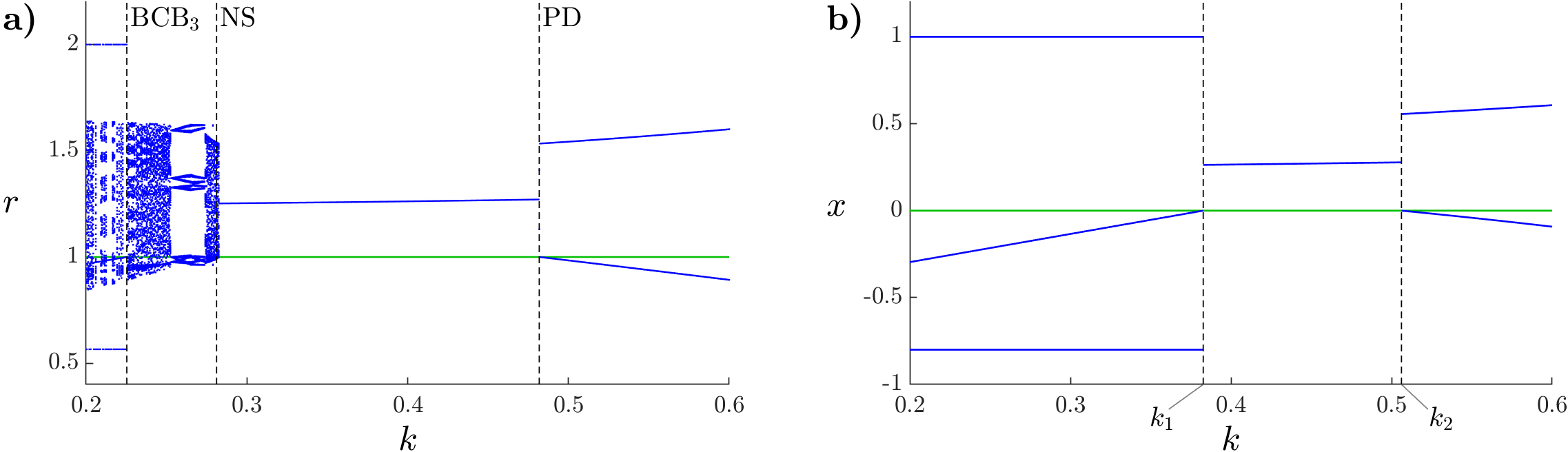}
\caption{
Panel (a) is a bifurcation diagram of the influenza model \eqref{eq:m10-RoHi19b} with $c = 0.9$ and $R_0 = 2$
(${\rm BCB}_3$:~BCB of a period-three solution;
NS:~Neimark-Sacker bifurcation;
PD:~period-doubling bifurcation).
Panel (b) is a bifurcation diagram of \eqref{eq:f} with \eqref{eq:oG} and $\mu = 1$
and corresponds to the influenza model in the limit $R_0 \to 1^+$.
The bifurcation values $k_1$ and $k_2$ are given by \eqref{eq:k1} and \eqref{eq:k2}.
\label{fig:zRRoHi19b_B}
} 
\end{center}
\end{figure}

Fig.~\ref{fig:zRRoHi19b_B}a is a bifurcation diagram
showing how the long-term dynamics of \eqref{eq:m10-RoHi19b} varies with $k$, using $c = 0.9$ and $R_0 = 2$.
For immediate values of $k$ the map has a stable fixed point
corresponding to the occurrence of an outbreak every year.
This fixed point loses stability in a supercritical period-doubling bifurcation at $k \approx 0.48182$,
shortly after which the period-doubled solution undergoes a BCB.
Indeed these two bifurcations occur so close together that on the scale of Fig.~\ref{fig:zRRoHi19b_B}a
it seems like the attractor changes discontinuously in a single bifurcation.

Beyond the BCB the period-doubled solution corresponds to outbreaks occurring every other year.
For relatively small values of $k$ the stability of the fixed point is lost in a Neimark-Sacker bifurcation.
There is also a period-adding structure and a stable period-three solution corresponding to the occurrence of an outbreak once every three years.

The basic skeleton of the bifurcation diagram can be explained from our above results for \eqref{eq:f}.
With $R_0 < 1$ the point $(S,T) = (0,1)$ is a stable fixed point
of \eqref{eq:m10-RoHi19b} (corresponding to the occurrence of no outbreaks).
At $R_0 = 1$ this fixed point undergoes a BCB by colliding with the switching manifold $r(S,T) = 1$.
The map \eqref{eq:m10-RoHi19b} is continuous across this manifold,
indeed the piecewise-linear approximation to the map centred at the bifurcation
can be converted to the normal form \eqref{eq:bcnf} via an affine change of coordinates.
By differentiating \eqref{eq:m10-RoHi19b}, we obtain
\begin{equation}
\tau_L = 0, \quad
\delta_L = 0, \quad
\tau_R = -2 c, \quad
\delta_R = 2 (1-k) c^2,
\label{eq:oG}
\end{equation}
where the left piece of \eqref{eq:bcnf} corresponds to the $r(S,T) < 1$ piece of \eqref{eq:m10-RoHi19b},
and the right piece of \eqref{eq:bcnf} corresponds to the $r(S,T) > 1$ piece of \eqref{eq:m10-RoHi19b}.
Notice $\delta_L = 0$, because the $r(S,T) < 1$ piece of \eqref{eq:m10-RoHi19b} has degenerate range,
hence the normal form is an instance of \eqref{eq:f}.

We now discuss the dynamics of \eqref{eq:f} with \eqref{eq:oG}.
With $\mu < 0$ the fixed point $(x^L,y^L)$ is admissible and asymptotically stable because $|\tau_L| < 1$.
With $\mu > 0$ the dynamics are shown in Fig.~\ref{fig:zRRoHi19b_B}b.
Notice from \eqref{eq:oG} that an increase in the value of $k$ corresponds to a decrease in the value of $\delta_R$.
So in terms of the $(\tau_R,\delta_R)$-plane (several of which were shown in \S\ref{sec:pos}, although none used $\tau_L = 0$),
moving left to right across Fig.~\ref{fig:zRRoHi19b_B}b corresponds to moving downwards along a vertical line in the $(\tau_R,\delta_R)$-plane.
The interval of $k$-values that give a stable fixed point
are where this line intersects the purple triangle $|\tau_R| - 1 < \delta_R < 1$.
The left endpoint of this interval is where $\delta_R = 1$, so by \eqref{eq:oG}
\begin{equation}
k_1 = 1 - \frac{1}{2 c^2}.
\label{eq:k1}
\end{equation}
To the left of this value \eqref{eq:f} has a stable $LLR$-cycle.
The right endpoint of the interval is where $\delta_R = -\tau_R - 1$, so by \eqref{eq:oG}
\begin{equation}
k_2 = 1 - \frac{1}{c} + \frac{1}{2 c^2},
\label{eq:k2}
\end{equation}
and to the right of this value \eqref{eq:f} has a stable $LR$-cycle.

The dynamics shown in Fig.~\ref{fig:zRRoHi19b_B}b correspond to those of the influenza model with a value of $R_0$ just greater than $1$.
We observe these dynamics have the same features as Fig.~\ref{fig:zRRoHi19b_B}a, which uses $R_0 = 2$, except without period-adding.
In this way our understanding of \eqref{eq:f}
goes a long way to explaining the dynamics of the influenza model,
even at parameter values relatively far beyond the BCB at $R_0 = 1$ where \eqref{eq:f} was derived.

\section{Discussion}
\label{sec:conc}

In this paper we have provided a broad summary of the dynamics of \eqref{eq:f}.
The intention is for this work to be a resource to be used
to obtain basic understandings of BCBs in mathematical models.
The methodology for doing this was illustrated in \S\ref{sec:appl}.
For an ODE model the first step is formulate a return map, e.g.~a Poincar\'e map,
that captures the oscillatory dynamics of interest.
For a given BCB, one then evaluates the
derivative (Jacobian matrix) of the left and right pieces of the Poincar\'e map at the bifurcation point.
Except for the simplest models, these derivatives need to be evaluated numerically, e.g.~with finite differences.
If the map is two-dimensional, one evaluates the trace and determinants of these matrices
to obtain values for $\tau_L$, $\delta_L$, $\tau_R$, and $\delta_R$
(for higher dimensional maps refer to \cite{Di03,Si16}).
If one of the determinants of zero, we can enforce $\delta_L = 0$, see Remark \ref{rm:LRswap}, thus giving \eqref{eq:f}.

The different signs of $\mu$ correspond to different sides of the BCB.
For $\mu < 0$ the map \eqref{eq:f} has a stable fixed point if $|\tau_L| < 1$,
while if it has a stable periodic solution then this solution has exactly one point in the left-half plane.
For $\mu > 0$ attracting invariant circles are possible.
On these the dynamics is either mode-locked or quasi-periodic.
Figs.~\ref{fig:zRbifSetA}, \ref{fig:zRbifSetB}, \ref{fig:zRbifSetC}, and \ref{fig:zRbifSetD}
summarise the extent to which the dynamics of \eqref{eq:f} with $\mu > 0$ varies with $\tau_L$, $\tau_R$, and $\delta_R$.

Overall the range of the dynamics of \eqref{eq:f}
is fairly comparable to \eqref{eq:bcnf} with $\delta_L \ne 0$ and $\delta_R \ne 0$, \cite{BaGr99,GhMc23,SiMe08b,SuGa08}.
Robust chaos occurs widely, and there are many areas of parameter space where multiple attractors coexist.
It is interesting that many bifurcation curves of \eqref{eq:f}
admit closed-form expressions, while if $\delta_L \ne 0$ and $\delta_R \ne 0$ such curves are analytically intractable.
Most notably, shrinking point curves of \eqref{eq:f} are zeros of polynomial functions of $\tau_L$, $\tau_R$, and $\delta_R$,
whereas if $\delta_L \ne 0$ and $\delta_R \ne 0$ they typically have many (possibly a dense set) of points where they are non-differentiable.

In three places the two-parameter bifurcation diagrams show structures for which no general theory is known.
These are illustrated in Figs.~\ref{fig:zRbifSetCzoom}b, \ref{fig:zRbifSetCzoomAlt}, and \ref{fig:zRbifSetDzoom}.
We anticipate that such theory is achievable by analysing an induced map $\phi$.
This map can be defined by $x' = \phi(x)$, where $(x',0)$ is the next point in the forward orbit of $(x,0)$ under \eqref{eq:f}
that lies on the $x$-axis.
Since $\delta_L = 0$, any point in the left-half plane maps to the $x$-axis,
so the induced map captures all dynamics of \eqref{eq:f} except those constrained in the right-half plane (which are trivial because $f_R$ is affine).
Since $\phi$ is one-dimensional it is more amenable to an exact analysis.
It remains for future work to use this map to explain the unexplained bifurcation structures.

\section*{Acknowledgements}

This work was supported by Marsden Fund contract MAU2209 managed by Royal Society Te Ap\={a}rangi.

\end{document}